\documentclass[11pt]{article}

\usepackage{amsmath,amssymb}
\usepackage{graphicx, citesort}

\newcommand{\suptf}{\Gamma}     
\newcommand{\Rtf}{R_\suptf}             
\newcommand{\supt}{T}                   
\newcommand{\supf}{\Omega}              
\newcommand{\Rt}{R_\supt}               
\newcommand{\Rf}{R_\supf}               
\newcommand{\idt}{t}                            
\newcommand{\idf}{\omega}               
\newcommand{\idtf}{\gamma}              
\newcommand{\Auxt}{\mathcal{H}_\supt}
\newcommand{\ids}{s}                              
\newcommand{\dv}{S}               
\newcommand{\FT}{F}                             

\newcommand{\ich}{\mu}                  

\newcommand{\at}{\alpha_1}              
\newcommand{\af}{\alpha_2}              
\newcommand{\atf}{\alpha}
\newcommand{\Mt}{\Phi_1}                        
\newcommand{\Mf}{\Phi_2}                        
\newcommand{\Mtf}{\Phi}
\newcommand{\supat}{\Gamma_1}
\newcommand{\supaf}{\Gamma_2}
\newcommand{\Rat}{R_{\supat}}
\newcommand{\Raf}{R_{\supaf}}

\newcommand{\supatf}{\Gamma}
\newcommand{\vt}{\phi}                  
\newcommand{\vf}{\psi}                  
\newcommand{\indt}{I_1}
\newcommand{\indf}{I_2}

\newcommand{\idtp}{{\idt^\prime}}
\newcommand{\idfp}{{\idf^\prime}}
\newcommand{\alphap}{\alpha^\prime}
\newcommand{\suptfp}{{\suptf^\prime}}
\newcommand{\Rtfp}{R_{\suptfp}}
\newcommand{\dvp}{{\dv^\prime}}
\newcommand{\<}{\langle}
\renewcommand{\>}{\rangle}
\newcommand{\bpm}{\begin{pmatrix}}
\newcommand{\epm}{\end{pmatrix}}
\def\E{{\hbox{\bf E}}}
\def\P{{\hbox{\bf P}}}

\newcommand{\C}{\mathbb{C}}
\newcommand{\R}{\mathbb{R}}

\newcommand{\sgn}{\operatorname{sgn}}
\newcommand{\supp}{\operatorname{supp}}
\renewcommand{\Re}{\operatorname{Re}}
\renewcommand{\Im}{\operatorname{Im}}
\newcommand{\nullsp}{\operatorname{Null}}
\newcommand{\dmn}{\operatorname{dim}}
\newcommand{\trace}{\operatorname{Tr}}
\newcommand{\range}{\operatorname{Range}}
\newcommand{\re}{{\sf r}}
\newcommand{\im}{{\sf i}}

\newtheorem{theorem}{Theorem}[section]
\newtheorem{lemma}{Lemma}[section]
\newtheorem{corollary}{Corollary}[section]
\def \endprf{\hfill {\vrule height6pt width6pt depth0pt}\medskip}
\newenvironment{proof}{\noindent {\bf Proof} }{\endprf\par}

\numberwithin{equation}{section}

\begin{document}

\title{\vspace{-5mm} Quantitative Robust Uncertainty Principles and \\ 
Optimally Sparse Decompositions}

\author{Emmanuel J. Cand\`es and Justin  Romberg\\[2mm]
Applied and Computational Mathematics, Caltech, Pasadena, CA  91125}


\maketitle

\begin{abstract}
  
  In this paper, we develop a robust {\em uncertainty principle} for
  finite signals in $\C^N$ which states that for almost all choices
  $\supt,\supf\subset\{0,\ldots,N-1\}$ such that
\[
|\supt| + |\supf| \asymp (\log N)^{-1/2} \cdot N, 
\]
there is no signal $f$ supported on $\supt$ whose discrete Fourier
transform $\hat{f}$ is supported on $\supf$.  In fact, we can make the
above uncertainty principle {\em quantitative} in the sense that if
$f$ is supported on $\supt$, then only a small percentage of the
energy (less than half, say) of $\hat{f}$ is concentrated on $\supf$.

As an application of this robust uncertainty
principle (QRUP), we consider the problem of decomposing a signal into
a sparse superposition of spikes and complex sinusoids
\[
f(\ids) = \sum_{\idt\in\supt}\at(\idt)\delta(s-\idt) ~+~
\sum_{\idf\in\supf}\af(\idf)e^{i2\pi \idf \ids/N}/\sqrt{N}. 
\]
We show that if a generic signal $f$ has a decomposition $(\at,\af)$
using spike and frequency locations in $\supt$ and $\supf$
respectively, and  obeying
\[
|\supt| + |\supf| \leq \mathrm{Const}\cdot (\log N)^{-1/2}\cdot N, 
\] 
then $(\at, \af)$ is the {\em unique sparsest} possible decomposition
(all other decompositions have more non-zero terms).  In addition, if
\[
|\supt| + |\supf| \leq \mathrm{Const}\cdot (\log N)^{-1}\cdot N, 
\]
then the sparsest $(\at,\af)$ can be found by solving a convex
optimization problem. 

Underlying our results is a new probabilistic approach which insists
on finding the correct uncertainty relation or the optimally sparse
solution for nearly all subsets but not necessarily all of them, and
allows to considerably sharpen previously known results
\cite{DonohoStark,DonohoHuo}. In fact, we show that the fraction of
sets $(T, \Omega)$ for which the above properties do not hold can be
upper bounded by quantities like $N^{-\alpha}$ for large values of
$\alpha$.

The QRUP (and the application to sparse approximation) can be extended
to general pairs of orthogonal bases $\Mt,\Mf$ of $\C^N$.  For almost
all choices $\supat,\supaf\subset\{0,\ldots,N-1\}$ obeying
\[
|\supat| + |\supaf| \asymp \ich(\Mt,\Mf)^{-2} \cdot (\log N)^{-m},
\]
there is no signal $f$ such that $\Mt f$ is supported on $\supat$ and
$\Mf f$ is supported on $\supaf$ where $\ich(\Mt,\Mf)$ is the {\em
  mutual incoherence} between $\Mt$ and $\Mf$.
\end{abstract}


{\bf Acknowledgments.} E.~C.~is partially supported by National
Science Foundation grants DMS 01-40698 (FRG) and ACI-0204932 (ITR),
and by an Alfred P.  Sloan Fellowship.  J.~R.~is supported by those
same National Science Foundation grants. E.~C.~would like to thank
Terence Tao for helpful conversations related to this project. These
results were presented at the International Conference on
Computational Harmonic Analysis, Nashville, Tennessee, May 2004.

\section{Introduction}

\subsection{Uncertainty principles}
\label{sec:introup}

The classical Weyl-Heisenberg uncertainty principle states that a
continuous-time signal cannot be simultaneously well-localized in both
time and frequency.  Loosely speaking, this principle says that if
most of the energy of a signal $f$ is concentrated near a
time-interval of length $\Delta t$ and most of its energy in the
frequency domain is concentrated near an interval of length $\Delta
\omega$, then
\[
\Delta t \cdot \Delta \omega \ge 1. 
\]
This principle is one of the major intellectual achievements of the
20th century and since then, much work has been concerned with
extending such uncertainty relations to other setups, namely, by
investigating to what extent it is possible to concentrate a function
$f$ and its Fourier transform $\hat f$, relaxing the assumption that
$f$ and $\hat f$ be concentrated near intervals as in the work of
Landau, Pollack and Slepian \cite{prolateI,prolateII,prolateIII}, or
by considering signals supported on a discrete set
\cite{prolateV,DonohoStark}.

Because our paper is concerned with finite signals, we now turn our
attention to ``discrete uncertainty relations'' and begin by recalling
the definition of the discrete Fourier transform
\begin{equation}
  \label{eq:dft}
  \hat f(\idf) = \frac{1}{\sqrt{N}} \sum_{t =0}^{N-1} f(t)  
e^{-i 2\pi \idf t/N},
\end{equation}
where the frequency index $\idf$ ranges over the set $\{0, 1, \ldots,
N - 1\}$. For signals of length $N$, \cite{DonohoStark} introduced a
sharp uncertainty principle which simply states that the supports of a
signal $f$ in the time and frequency domains must obey
\begin{equation}
\label{eq:discreteup}
|\supp f| + |\supp \hat{f}| \geq 2\sqrt{N}.
\end{equation}
We emphasize that there are no other restriction on the organization
of the supports of $f$ and $\hat f$ other than the size constraint
\eqref{eq:discreteup}. \cite{DonohoStark} also observed that the
uncertainty relation \eqref{eq:discreteup} is tight in the sense that
equality is achieved for certain special signals.  For example,
consider as in \cite{DonohoStark,DonohoHuo} the {\em Dirac comb}
signal: we suppose that the sample size $N$ is a perfect square and
let $f$ be equal to 1 at multiples of $\sqrt{N}$ and 0
everywhere else
\begin{equation}
  \label{eq:dirac}
  f(t) = \begin{cases}
1, & t=m\sqrt{N},~~ m=0,1,\ldots,\sqrt{N}-1 \\
0, & \text{elsewhere}. 
\end{cases} 
\end{equation}
Remarkably, the Dirac comb is invariant through the Fourier transform,
i.e. $\hat{f} = f$, and therefore, $|\supp f |+|\supp \hat f| =
2\sqrt{N}$.  In other words, \eqref{eq:discreteup} holds with
equality.

In recent years, uncertainty relations have become very popular, in
part because they help explaining some miraculous properties of
$\ell_1$-minimization procedures as we will see below, and researchers
have naturally developed similar uncertainty relations between pairs
of bases other than the canonical basis and its conjugate. We single
out the work of Elad and Bruckstein \cite{EladBruckstein} which
introduces a generalized uncertainty principle for pairs $\Mt,\Mf$ of
orthonormal bases.  Define the {\em mutual incoherence}
\cite{DonohoHuo,GribonvalNielsen,EladBruckstein} between $\Phi_1$ and
$\Phi_2$ as
\begin{equation}
\label{eq:mich}
\ich(\Mt,\Mf) = 
\max_{\vt\in\Mt, \vf\in\Mf} |\<\vt,\vf\>|; 
\end{equation}
then if $\at$ is the (unique) representation of $f$ in basis $\Mt$
with $\supat=\supp\at$, and $\af$ is the representation in $\Mf$, the
supports must obey
\begin{equation}
\label{eq:gendiscreteup}
|\supat| + |\supaf| \geq \frac{2}{\ich(\Mt,\Mf)}. 
\end{equation} 
Note that the mutual incoherence $\mu$ always obeys $1/\sqrt{N} \le
\mu \le 1$ and measures how the two bases look alike. The smaller the
incoherence, the stronger the uncertainty relation.   To see how this generalizes the discrete uncertainty principle, observe that in the case where $\Mt$ is the canonical or spike basis and $\Mf$ is the Fourier basis, $\ich = 1/\sqrt{N}$ (maximal incoherence) and \eqref{eq:gendiscreteup} is, of course, \eqref{eq:discreteup}.

\subsection{The tightness of the uncertainty relation is fragile}

It is true that there exist signals that saturate the uncertainty
relations but such signals are very special and are hardly
representative of ``generic'' or ``most'' signals. Consider the Dirac
comb for instance; here the locations and heights of the $\sqrt{N}$
spikes in the time domain carefully conspire to create an inordinate
number of cancellations in the frequency domain.  This will not be the
case for sparsely supported signals in general.  Simple numerical
experiments confirm that signals with the same support as the Dirac
comb but with different spike amplitudes almost always have Fourier
transforms that are nonzero everywhere.  Indeed, constructing pathological examples other than the Dirac comb requires mathematical wizardry.

Moreover, if the signal length $N$ is prime (making signals like the
Dirac comb impossible to construct), the discrete uncertainty
principle is sharpened to \cite{tao:uncertainty}
\begin{equation}
\label{eq:primeup}
|\supp f| + |\supp \hat{f}| > N,
\end{equation}
which validates our intuition about the exceptionality of signals
such as the Dirac comb. 

\subsection{Robust uncertainty principles}

Excluding these exceedingly rare and exceptional pairs $T: = \supp f, \Omega: = \supp\hat f$, how tight is the uncertainty relation? That is, given two
sets $T$ and $\Omega$, how large need $|T| + |\Omega|$ be so that it
is possible to construct a signal whose time and frequency supports
are $T$ and $\Omega$ respectively? In this paper, we introduce a {\em
  robust} uncertainty principle (for general $N$) which illustrates
that for ``most'' sets $\supt,\supf$, \eqref{eq:primeup} is closer to
the truth than \eqref{eq:discreteup}.  Suppose that we choose 
$(\supt,\supf)$ at random from all pairs obeying
\begin{equation*}
  |\supt|+|\supf| \le  \frac{N}{\sqrt{(\beta + 1) \log N}}.
\end{equation*}
Then with overwhelming high probability---in fact, exceeding $1 -
O(N^{-\beta \rho})$ for some positive constant $\rho$ (we shall give
explicit values)---we will be unable to find a signal in $\C^N$
supported on $\supt$ in the time domain and $\supf$ in the frequency
domain. In other words, remove a negligible fraction of sets and
\begin{equation}
\label{eq:rup}
  |\supp f|+|\supp \hat f| >  \frac{N}{\sqrt{(\beta + 1) \log N}}, 
\end{equation}
holds, not \eqref{eq:discreteup}.

Our uncertainty principle is not only robust in the sense that it
holds for most sets, it is also {\em quantitative}. Consider a random
pair $(T, \Omega)$ as before and put $1_\Omega$ to be the indicator
function of the set $\Omega$. Then with essentially the same
probability as above, we have
\begin{equation}
  \label{eq:qrup}
\|\hat f \cdot 1_\Omega\|^2 \le \|\hat f\|^2/2,  
\end{equation}
say, for all functions $f$ supported on $T$. By symmetry, the same
inequality holds by exchanging the role of $T$ and $\Omega$,
\[
\|f \cdot 1_T\|^2 \le \|f\|^2/2,
\]
for all functions $\hat f$ supported on $\Omega$.  Moreover, 
as with the discrete uncertainty principle, the QRUP can be extended to arbitrary pairs of bases.

\subsection{Significance of uncertainty principles}

In the last three years or so, there has been a series of papers
starting with \cite{DonohoHuo} establishing a link between discrete
uncertainty principles and sparse approximation
\cite{DonohoHuo,GribonvalNielsen,DonohoElad,Tropp03}. In this field, the goal
is to separate a signal $f \in \C^N$ into two (or more) components,
each representing contributions from different phenomena. The idea is
as follows: suppose we have two (or possibly many more) orthonormal
bases $\Mt,\Mf$; we search among all the decompositions $(\at,\af)$ of
the signal $f$
\[
f = \bpm \Mt & \Mf \epm \bpm \at \\ \af \epm := \Mtf\atf
\]
for the shortest one
\begin{equation}
  \label{eq:P_0}
(P_0) \quad\quad \min_{\atf} \|\atf\|_{\ell_0},\quad \Mtf\atf = f, 
\end{equation}
where $\|\atf\|_{\ell_0}$ is simply the size of the support of
$\alpha$, $\|\atf\|_{\ell_0} := |\{\idtf,\ \atf(\idtf) \neq 0\}|$.

The discrete uncertainty principles \eqref{eq:discreteup} and
\eqref{eq:gendiscreteup} are useful in the sense that they tell us
when $(P_0)$ has a unique solution.  When $\Mtf$ is the time-frequency
dictionary, it is possible to show that if a signal $f$ has a
decomposition $f = \Mtf\atf$ consisting of spikes on subdomain $\supt$
and frequencies on $\supf$, and
\begin{equation}
\label{eq:dhl0}
|\supt| + |\supf| < \sqrt{N}, 
\end{equation}
then $\atf$ is the unique minimizer of $(P_0)$ \cite{DonohoHuo}. In a
nutshell, the reason is that if $\Phi(\alpha_0 + \delta_0)$ were
another decomposition, $\delta_0$ would obey $\Phi \delta_0 = 0$ which
says that $\delta_0$ would be of the form $\delta_0 = (\delta, - \hat
\delta)$. Now \eqref{eq:discreteup} implies that $\delta_0$ would have
at least $2\sqrt{N}$ nonzero entries which in turn would give
$\|\alpha_0 + \delta_0\|_{\ell_0} \ge \sqrt{N}$ for all $\alpha_0$
obeying $\|\alpha_0\|_{\ell_0} < \sqrt{N}$---thereby proving the
claim.  Note that again the condition \eqref{eq:dhl0} is sharp because
of the extremal signal \eqref{eq:dirac}. Indeed, the Dirac comb may be
expressed as a superposition of $\sqrt{N}$ terms in the time or in the
frequency domain; for this special signal, $(P_0)$ does not have a
unique solution.

In \cite{EladBruckstein}, the same line of reasoning is followed for
general pairs of orthogonal bases, and $\ell_0$-uniqueness is
guaranteed when
\begin{equation}
\label{eq:ebl0}
|\supat| + |\supaf| < \frac{1}{\ich(\Mt,\Mf)}.
\end{equation}

Unfortunately, as far as finding the sparsest decomposition, solving
$(P_0)$ directly is computationally infeasible because of the highly
non-convex nature of the $\|\cdot\|_{\ell_0}$ norm.  To the best of
our knowledge, finding the minimizer obeying the constraints would
require searching over all possible {\em subsets} of columns of
$\Mtf$, an algorithm that is combinatorial in nature and has
exponential complexity.  Instead of solving $(P_0)$, we consider a
similar program in the $\ell_1$ norm which goes by the name of {\em
  Basis Pursuit} \cite{BP}:
\begin{equation}
\label{eq:(P_1)}
(P_1)\quad \quad  \min_{\atf} 
\|\atf\|_{\ell_1}, \quad \Mtf\atf = f.
\end{equation}
Unlike the $\ell_0$ norm, the $\ell_1$ norm is convex.  As a result,
$(P_1)$ can be solved efficiently using standard ``off the shelf''
optimization algorithms.  The $\ell_1$-norm can also be viewed as a
``sparsity norm'' which among the vectors that meet the constraints,
will favor those with a few large coefficients and many small
coefficients over those where the coefficient magnitudes are
approximately equal \cite{BP}.

A beautiful result in \cite{DonohoHuo} actually shows that if $f$ has
a sparse decomposition $\atf$ supported on $\supatf$ with
\begin{equation}
\label{eq:dhgenl1}
|\supatf| < \frac{1}{2}(1+\ich^{-1}), 
\end{equation}
then the minimizer of $(P_1)$ is unique and is equal to the minimizer
of $(P_0)$ (\cite{EladBruckstein} improves the constant in
\eqref{eq:dhgenl1} from $1/2$ to $\approx .9142$).  In these
situations, we can replace the highly non-convex program $(P_0)$ with
the much tamer (and convex) $(P_1)$.

We now review a few applications of these types of ideas.
\begin{itemize}
\item {\em Geometric Separation.} Suppose we have a dataset and one
  wishes to separate point-like structures, from filamentary
  (edge-like) structures, from sheet-like structures. In 2 dimensions,
  for example, we might imagine synthesizing a signal as a
  superposition of wavelets and curvelets which are ideally adapted to
  represent point-like and curve-like structures respectively.
  Delicate space/orientation uncertainty principles show that the
  minimum $\ell_1$-norm decomposition in this combined dictionary
  automatically separates point and curve-singularities; the wavelet
  component in the decomposition \eqref{eq:(P_1)} accurately captures
  all the pointwise singularities, while the curvelet component
  captures all the edge curves. We refer to \cite{DonohoGeomSep} for
  theoretical developments and to \cite{StarckAstro} for numerical
  experiments.
  
\item {\em Texture-edges separation} Suppose now that we have an image
  we wish to decompose as a sum of a cartoon-like geometric part plus a
  texture part.  The idea again is to use curvelets to represent the
  geometric part of the image and local cosines to represent the
  texture part. These ideas have recently been tested in practical settings, 
  with spectacular success \cite{EdgesTextures} (see also
  \cite{MeyerAverbuchCoifman} for earlier and related ideas).
\end{itemize}

In a different direction, the QRUP is also implicit in some
of our own work on the exact reconstruction of sparse signals from
vastly undersampled frequency information \cite{CRT}.  Here, we wish
to reconstruct a signal $f \in \C^N$ from the data of only $|\Omega|$
random frequency samples. The surprising result is that although most
of the information is missing, one can still reconstruct $f$ {\em
  exactly} provided that $f$ is sparse.  Suppose $|\Omega|$ obeys the
oversampling relation
\[
|\Omega| \asymp |T| \cdot \log N
\]
with $T:= \supp f$. Then with overwhelming probability, the object $f$
(digital image, signal, and so on)  is the exact and unique solution of the
convex program that searches, among all signals that are consistent
with the data, for that with minimum $\ell_1$-norm. 
We will draw on the the tools developed in the earlier work, making the QRUP {\em explicit} and applying it to the problem of searching for sparse decompositions.

\subsection{Innovations}

Nearly all the existing literature on uncertainty relations and its
consequences focuses on worst case scenarios, compare
\eqref{eq:discreteup} and \eqref{eq:dhgenl1}. What is new here is the
development of probabilistic models which show that the performance of
Basis Pursuit in an overwhelmingly large majority of situations is
actually very different than that predicted by the overly
``pessimistic'' bounds \eqref{eq:dhgenl1}.  For the time-frequency
dictionary, we will see that if a representation $\atf$ (with spike
locations $\supt$ and sinusoidal frequencies $\supf$) of a signal $f$
exists with
\[
|\supt| + |\supf| \asymp N/\sqrt{\log N},
\]
then $\atf$ is the sparsest representation of $f$ almost all of the
time.  If in addition, $\supt$ and $\supf$ satisfy
\begin{equation}
  \label{eq:sparserup}
  |\supt|+|\supf| \asymp N/\log N,
\end{equation}
then $\atf$ can be recovered by solving the convex program $(P_1)$.
In fact, numerical simulations reported in section \ref{sec:numerical}
suggest that \eqref{eq:sparserup} is far closer to the empirical
behavior than \eqref{eq:dhgenl1}, see also \cite{DonohoStark}. We show
that similar results also hold for general pairs of bases $\Mt,\Mf$.

As discussed earlier, there is by now a well-established machinery
that allows turning uncertainty relations into statements about the
ability to find sparse decompositions. We would like to point out that
our results \eqref{eq:sparserup} are not an automatic consequence of
the uncertainty relation \eqref{eq:rup} together with these existing
ideas.  Instead, our analysis relies on the study of eigenvalues of
random matrices which, of course, is completely new.

\subsection{Organization of the paper} 

In Section~\ref{sec:pmodel} we develop a probability model that shall
be used throughout the paper to formulate our results. In Section
\ref{sec:qrup}, we will establish uncertainty relations such as
\eqref{eq:qrup}. Sections \ref{sec:bpIF} and \ref{sec:BPgeneral} will
prove uniqueness and equality of the $(P_0)$ and $(P_1)$ programs.  In
the case where the basis pair $(\Phi_1, \Phi_2)$ is the time-frequency
dictionary (Section \ref{sec:bpIF}), we will be very careful in
calculating the constants appearing in the bounds.  We will be
somewhat less precise in the general case (Section
\ref{sec:BPgeneral}), and will forgo explicit calculation of
constants. We report on numerical experiments in Section
\ref{sec:numerical} and close the paper with a short discussion
(Section \ref{sec:discussion}).

\section{A Probability Model for $\supat,\supaf$}
\label{sec:pmodel}

To state our results precisely, we first need to specify a
probabilistic model. We let $\indt$ and $\indf$ be two independent
Bernoulli sequences with parameters $p_T$ and $p_\Omega$ respectively
\begin{eqnarray*}
\indt(\idt)  = 1 & ~~\text{with probability}~~ p_T\\
\indf(\idf)  =  1 & ~~\text{with probability}~~ p_\Omega
\end{eqnarray*}
where $\idt,\idf=0,\ldots,N-1$, and define the support sets for the
spikes and sinusoids (and in general for the bases $\Mt$ and $\Mf$)
as
\begin{equation}
\label{eq:random}
\supt = \{\idt~~\text{s.t.}~~\indt(\idt) = 1 \}, \qquad 
\supf = \{\idf~~\text{s.t.}~~\indf(\idf) = 1\}. 
\end{equation}
If both $p_T$ and $p_\Omega$ are not too small, an application of the 
standard large deviations inequality shows us that our model is 
approximately equivalent to sampling $\E |T| = p_T \cdot N$ spike
locations and $\E |\Omega| = p_\Omega \cdot N$ frequency locations
uniformly at random. 

As we will see in the next section, the robust uncertainty principle
holds---with overwhelming probability---over sets $\supt$ and $\supf$
randomly sampled as above. Our estimates are quantitative and
introduce sufficient conditions so that the probability of ``failure''
be arbitrarily small, i.e. less than $O(N^{-\beta})$ for some
arbitrary $\beta > 0$. As a consequence, we will always assume that
\begin{equation}
  \label{eq:betalogN}
  \min(\E |T|, \E |\Omega|) \ge 4 (\beta+1) \cdot \log N
\end{equation}
as otherwise, one would have to consider situations in which $T$ or
$\Omega$ (or both) are empty sets---a situation of rather limited
interest. We also note that for $p_T$ and $p_\Omega$ as above, we have 
\begin{equation}
  \label{eq:2T}
  \P(|T| > 2 p_T \cdot N) \le N^{-\beta},
\end{equation}
as this follows from the well-known large deviation bound \cite{MassartSharp}
$$\P(|T| > \E |T| + t) \le \exp\left(-\frac{t^2}{2 \E |T| +
    2t/3}\right).$$

Further, to establish sparse approximation bounds (section
\ref{sec:bpIF}), we will also introduce a probability model on the
``active'' coefficients.  Given a pair $(\supt, \supf)$, we sample the
coefficient vector $\{\atf(\idtf),~\idtf\in\supatf\}$ from a
distribution with identically and independently distributed
coordinates; we also impose that each $\atf(\idtf)$ be drawn from a
continuous probability distribution that is {\em circularly symmetric}
in the complex plane; that is, the phase of $\atf(\idtf)$ is uniformly
distributed on $[0,2\pi)$.


\section{Quantitative Robust Uncertainty Principles}
\label{sec:qrup}

Equipped with the probability model \eqref{eq:random}, we now
introduce our uncertainty relations. To state our result, we make use
of the standard notation $o(1)$ to indicate a numerical term
tending to $0$ as $N$ goes to infinity.
\begin{theorem}
\label{th:qrupIF}
Assume the parameters in the model \eqref{eq:random} obey
\begin{equation}
\label{eq:tauqrup}
2\sqrt{\E |T| \cdot \E |\Omega|}  \le  \E |T| + \E |\Omega| \le 
\frac{N}{\sqrt{(\beta+1)\log N}} \, (\rho_0/2  + 
o(1)), \qquad \rho_0 = .7614
\end{equation}
(we will assume throughout the paper that $\beta \ge 1$ and $N \ge
512$) and let $(\supt,\supf)$ be a randomly sampled support pair.
Then with probability at least $1-O(\log N \cdot N^{-\beta})$; {\em
  every} signal $f$ supported on $\supt$ in the time domain has most
of its energy in the frequency domain outside of $\supf$
\[
\|\hat f \cdot 1_\Omega\|^2 \leq \frac{\|\hat f\|^2}{2};
\]
and likewise, {\em every} signal $f$ supported on $\supf$ in the
frequency domain has most of its energy in the time domain outside of
$\supt$
\[
\|\hat f \cdot 1_T\|^2 \leq \frac{\|f\|^2}{2}.
\]
As a result, it is impossible to find a signal $f$ supported on
$\supt$ whose discrete Fourier transform $\hat{f}$ is supported on
$\supf$. For finite sample sizes $N$, we can select
the parameters in \eqref{eq:tauqrup} as
\[
\E |T| + \E |\Omega| \le \frac{.2660\, N}{\sqrt{(\beta+1)\log N}}.
\]
\end{theorem} 
To establish this result, we introduce (as in \cite{CRT}) the
$|\supt|\times |\supt|$ auxiliary matrix $\Auxt$
\begin{equation}
\label{eq:auxt}
\Auxt(\idt,\idt^\prime) = \begin{cases} 
0 & \idt = \idt^\prime \\
\sum_{\idf\in\supf} e^{i\idf (\idt-\idt^\prime)} & \idt\not=\idt^\prime
\end{cases}.
\end{equation}
The following lemma effectively says that the eigenvalues of $\Auxt$
are small compared to $N$.
\begin{lemma}
\label{lm:auxeigs}
Fix $q$ in $(0,1)$ and suppose that 
\[
p_T + p_\Omega \le \rho_0 \cdot \frac{q}{\sqrt{(\beta+1)\log N}},
\qquad \rho_0 = .7614.
\]
Then the the matrix $\Auxt$ obeys
\[
\P\left( \|\Auxt\| \geq q N \right) \leq (\beta + 1) \log N \cdot
N^{-\beta}.
\] 
\end{lemma}
\begin{proof}
  The Markov inequality gives \begin{equation}
\label{eq:markov}
\P(\|\Auxt\| \geq q N) \leq \frac{\E\|\Auxt^n\|^2}{q^{2n}N^{2n}}, 
\quad  \text{for all } n \geq 1. 
\end{equation}
Recall next that the Frobenius norm $\|\cdot\|_F$ dominates the
operator norm $\|\Auxt\|\leq\|\Auxt\|_F$. This fact allows to leverage
results from \cite{CRT} which derives bounds for the conditional
expectation $\E[\|\Auxt^n\|^2_F \, | \, T]$ (where the expectation is
over $\supf$ for a {\em fixed} $\supt$):
\[ 
\E[\|\Auxt^n\|^2_F \, | \,T] \leq
(2n)\left(\frac{(1+\sqrt{5})^2}{2e(1-p_\Omega)}\right)^n n^n |\supt|^{n+1}
\, p_{\Omega}^n \, N^n.
\]
Our assumption about the size of $p_T + p_\Omega$ assures that
$p_\Omega < .12$ so that $(1+\sqrt{5})^2/2(1-p_\Omega) \le 6$, whence
\begin{equation}
\label{eq:crtbnd} 
\E[\|\Auxt^n\|^2_F \, | \,T] \leq 2n \, 
(6/e)^n n^n |\supt|^{n+1} \, p_\Omega^n \, N^n. 
\end{equation}
We will argue below that for $n \le (\beta + 1) \log N$ and $p_T$
obeying \eqref{eq:betalogN}
\begin{equation}
  \label{eq:chernoff}
  \E [|T|^{n+1}] \le 1.15^{n+1} \, [\E |T|]^{n+1} = 
1.15^{n+1} \, (p_T \, N)^{n+1}.  
\end{equation}
Since $p_T \le .25$, we established
\[
\E \|\Auxt^n\|^2_F \le (6 \times 1.15/e)^n \, n^{n+1} \cdot 
p_T^{n} \, p_\Omega^{n} \, N^{2n+1}.
\]
Observe now that together with $\sqrt{p_T \, p_\Omega} \le (p_T +
p_\Omega)/2$, this gives
\begin{equation}
\label{eq:pb1}
\P(\|\Auxt\| \geq q N)  \leq  
\left(\frac{p_T + p_\Omega}{\rho_0 \, q}\right)^{2n} \, e^{-n} \, n^{n+1}  N, 
\qquad \rho_0 = 1/\sqrt{6 \times 1.15} = .7614.  
\end{equation}
We now specialize \eqref{eq:pb1} and take $n= \lceil (\beta+1)\log N
\rceil$ where $\lceil x \rceil$ is the smallest integer greater or
equal to $x$. Then if $p_T + p_\Omega$ obeys \eqref{eq:tauqrup},
\begin{equation}
\P(\|\Auxt\| \geq q N)  \leq  [(\beta+1)\log N + 1] \cdot N^{-\beta}, 
\end{equation}
as claimed. 
\end{proof}

We now return to  \eqref{eq:chernoff} and write $|T|$ as
$$
|T| = \E |T| \cdot (1 + Y), \qquad Y = \frac{|T| - \E |T|}{\E |T|}. 
$$
Then 
$$
\E |T|^{n+1} = (\E |T|)^{n+1} \cdot \E (1 + Y)^{n+1} \le (\E
|T|)^{n+1} \cdot \E [\exp((n+1) Y)].
$$
Observe that $Y$ is a an affine function of a sum of independent
Bernoulli random variables. Standard calculations then give
\[
\E [\exp(n Y)] = e^{-n} \cdot \left(1 + \frac{n}{N}
  \frac{e^\lambda - 1}{\lambda} \right)^N, \qquad \lambda = n/(N p_T).
\]
Recall the assumption \eqref{eq:betalogN} which implies $\lambda \le
1/4$ which in turn gives $\lambda^{-1}(e^\lambda - 1) -
1 \le \log 1.15$. The claim follows. 

We would like to remark that \eqref{eq:chernoff} might be considerably
improved when $\E |T| = p_T \cdot N$ is much larger than $n$ since in
that case, the binomial will have enhanced concentration around its
mean. For example,
\[
\E |T|^{n+1} \le 2 \cdot [\E |T|]^{n+1}
\]
in the event where $n \le \rho \cdot \sqrt{p_T \, N}$ for some
positive constant $\rho$ that the above method allows to calculate
explicitely. This would of course lead to improved constants in
\eqref{eq:tauqrup} and in the statement of Lemma \ref{lm:auxeigs}. In
this paper, we shall not pursue all these refinements as not to
clutter the exposition.

\begin{proof}{\bf of Theorem~\ref{th:qrupIF}}
  Let $f\in\C^N$ be supported on $\supt$; as such, $\Rt^*\Rt f = f$,
  where $\Rt$ is the restriction operator to $\supt$.  Put $F_{\Omega
    T}= \Rf\FT\Rt^*$. We have
\begin{equation*}
\|\hat{f} \cdot 1_\Omega\|_2  =  
\|F_{\Omega T} f\|_2 \le \|F_{\Omega T}\| \cdot \|f\|_2, 
\end{equation*}
and since, $\|F_{\Omega T}\|^2 = \|F_{\Omega T}^*F_{\Omega T}\|$, it
will suffice to show that that with high probability, the largest
eigenvalue of $F_{\Omega T}^*F_{\Omega T}$ is less than $1/2$.

Using the definition of the auxiliary matrix in \eqref{eq:auxt}, it is
not hard to verify the identity $F_{\Omega T}^*F_{\Omega T} =
\frac{|\supf|}{N} I + \frac{1}{N}\Auxt$.  Suppose that $p_T + p_\Omega$ 
obeys 
the condition in Lemma~\ref{lm:auxeigs}; then 
except for a set of probability less than $O(\log N \cdot N^{-\beta})$,
\begin{equation*}
 \frac{|\supf|}{N} \le 
2 p_\Omega \le 2 \rho_0 \cdot \frac{q}{\sqrt{(\beta+1) \log N}}, 
\quad \text{ and } \quad \frac{1}{N}\|\Auxt\| \le q, 
\end{equation*}
and, therefore,
\begin{equation}
\label{eq:AtAnorm}
\|F_{\Omega T}^*F_{\Omega T}\| \leq q \cdot \left(1 + 
\frac{2 \rho_0}{\sqrt{(\beta+1) \log N}}\right) = q(1+o(1)).
\end{equation}
The theorem follows from taking $q = 1/2 + o(1)$. For the statement
about finite sample sizes, we observe that for $\beta \ge 1$ and $N
\ge 512$, $2/\sqrt{(\beta+1) \log N} \le .567$ and, therefore,
$\|F_{\Omega T}^*F_{\Omega T}\| \le 1/2$ provided that $q \le [2(1+
.567\rho_0)]^{-1}$.  This establishes the first part of the theorem.

By symmetry of the discrete Fourier transform, the claim about the
size of $\|f \cdot 1_T\|$ for $\hat f$ supported on a random set
$\Omega$ is proven exactly in the same way. This concludes the proof
of the theorem.
\end{proof}


\section{Robust UPs and Basis Pursuit: Spikes and Sinusoids} 
\label{sec:bpIF}

As in \cite{DonohoHuo,GribonvalNielsen,EladBruckstein}, our
uncertainty principles are directly applicable to finding sparse
approximations in redundant dictionaries.  In this section, we look
exclusively at the case of spikes and sinusoids. We will leverage
Theorem~\ref{th:qrupIF} in two different ways:
\begin{enumerate}
\item $\ell_0$-uniqueness: If $f\in\C^N$ has a decomposition $\atf$
  supported on $\supt\cup\supf$ with $|\supt|+|\supf|\asymp (\log
  N)^{-1/2} N$, then with high probability, $\atf$ is the sparsest
  representation of $f$.
\item Equivalence of $(P_0)$ and $(P_1)$: If $|\supt|+|\supf|\asymp
  (\log N)^{-1} N$, then $(P_1)$ recovers $\atf$ with overwhelmingly
  large probability.
\end{enumerate}

\subsection{$\ell_0$-uniqueness}
\label{sec:l0uniq}

To illustrate that it is possible to do much better than
\eqref{eq:dhl0}, we first consider the case in which $N$ is a prime
integer.  Tao \cite{tao:uncertainty} derived the following exact,
sharp discrete uncertainty principle.
\begin{lemma}
\label{lm:primeup}\cite{tao:uncertainty} Suppose that the sample size $N$ is a
prime integer. Then
\[
|\supp f| + |\supp \hat{f}| > N, \qquad \forall f \in \C^N. 
\]
\end{lemma} 

Using Lemma~\ref{lm:primeup}, a strong $\ell_0$-uniqueness result
immediately follows:
\begin{corollary}
\label{th:l0prime}
Let $\supt$ and $\supf$ be subsets of $\{0, \ldots, N-1\}$ for $N$
prime, and let $\alpha$ (with $\Phi\alpha = f$) be a vector supported
on $\suptf=\supt\cup\supf$ such that
\[
|\supt| + |\supf| \le N/2.
\]
Then the solution to $(P_0)$ is unique and is equal to $\alpha$.
Conversely, there exist distinct vectors $\alpha_0, \alpha_1$ obeying
$|\supp\alpha_0|, |\supp\alpha_1| \leq N/2 +1$ and $\Phi\alpha_0 =
\Phi \alpha_1$.
\end{corollary}
\begin{proof}
  As we have seen in the introduction, one direction is trivial. If
  $\alpha_0 + \delta_0$ is another decomposition, then $\delta_0$ is
  of the form $\delta_0: = (\delta, - \hat \delta)$. Lemma
  \ref{lm:primeup} gives $\|\delta_0\|_{\ell_0} > N$ and thus
  $\|\alpha_0 + \delta_0\|_{\ell_0} \ge \|\delta\|_{\ell_0} -
  \|\alpha\|_{\ell_0} > N/2$. Therefore, $\|\alpha_0 +
  \delta_0\|_{\ell_0} > \|\alpha\|_{\ell_0}$.
  
  For the converse, we know that since $\Phi$ has rank at most $N$, we
  can find $\delta \neq 0$ with $|\supp\delta| = N + 1$ such that
  $\Phi \delta = 0$. (Note that it is of course possible to construct
  such $\delta$'s for any support of size greater than $N$). Consider
  a partition of $\supp\delta = \suptf_0 \cup \suptf_1$ where
  $\suptf_0$ and $\suptf_1$ are two disjoint sets with $|\suptf_0| =
  N/2+1$ and $|\suptf_1| = N/2$, say.  The claim follows by taking
  $\alpha_0 = \delta|_{\suptf_0}$ and $\alpha_1 = -
  \delta|_{\suptf_1}$.
\end{proof}

A slightly weaker statement addresses arbitrary sample sizes.
\begin{theorem}
\label{th:qrupl0IF}
Let $f = \Phi \alpha$ be a signal with support set $\suptf =
\supt\cup\supf$ and coefficients $\atf$ sampled as in
Section~\ref{sec:pmodel}, and with parameters obeying
\eqref{eq:tauqrup}.  Then with probability at least $1-O(\log N \cdot
N^{-\beta})$, the solution to $(P_0)$ is unique and equal to $\atf$.
\end{theorem}
To prove Theorem \ref{th:qrupl0IF}, we shall need the following lemma:
\begin{lemma}
\label{lm:dimnull}
Suppose  $\supt$ and $\supf$ are fixed subsets of $\{0,\ldots,N-1\}$, put 
$\suptf
= \supt\cup\supf$, and let $\Phi_\suptf:=\Mtf\Rtf^*$ be the
$N\times(|\supt|+|\supf|)$ matrix $\Phi_\suptf = \bpm \Rt^* &
\FT^*\Rf^*\epm$. Then 
\[
|\suptf| < 2N \qquad \Rightarrow \qquad \dmn(\nullsp(\Phi_\suptf)) <
\frac{|\suptf|}{2}.
\]
\end{lemma}
\newcommand{\FOT}{{F_{\Omega T}}}
\newcommand{\FOTstar}{{F^*_{\Omega T}}}
\begin{proof}
Obviously, 
\[
\dmn(\nullsp(\Phi_\suptf)) = \dmn(\nullsp(\Phi^*_\suptf \Phi_\suptf)),
\]
and we then write the $|\suptf|\times|\suptf|$ matrix $\Phi^*_\suptf
\Phi_\suptf$ as
\[
\Phi_\suptf^* \Phi_\suptf = \bpm I & \FOTstar
\\ \FOT & I \epm
\]
with $\FOT$ the partial Fourier transform from $T$ to $\Omega$ $\FOT:=
\Rf\FT\Rt^*$.  The dimension of the nullspace of $\Phi_\suptf^*
\Phi_\suptf $ is simply the number of eigenvalues of $\Phi_\suptf^*
\Phi_\suptf $ that are zero.  Put 
\[
G := I - \Phi_\suptf^* \Phi_\suptf = \bpm 0 & \FOTstar \\ \FOT & 0 \epm, \quad
\text{so that} \quad G^*G = \bpm \FOTstar \FOT & 0 \\ 0 & \FOT \FOTstar \epm.
\]
Letting $\lambda_j(\cdot)$ denotes the $j$th largest eigenvalue of a
matrix, observe that $\lambda_j(\Phi_\suptf^* \Phi_\suptf ) = 1 -
\lambda_j(G)$, and since $G$ is symmetric
\begin{equation}
\label{eq:Gtr}
\trace(G^*G)  = 
\lambda_1^2(G) + \lambda_2^2(G) + \cdots + \lambda_{|\supt|+|\supf|}^2(G).
\end{equation}
We also have that $\trace(G^*G) = \trace(\FOTstar\FOT) +
\trace(\FOT\FOTstar)$, so the eigenvalues in \eqref{eq:Gtr} will
appear in duplicate,
\begin{equation}
\label{eq:Gtr2}
\lambda_1^2(G) + \lambda_2^2(G) + \cdots + \lambda_{|\supt|+|\supf|}^2(G) = 
2\cdot(\lambda_1^2(\FOTstar\FOT) + \cdots + 
\lambda_{|\supt|}^2(\FOTstar\FOT)).
\end{equation}
We calculate 
\[
(\FOTstar\FOT)_{\idt,\idtp} = \frac{1}{N}\sum_{\idf\in\supf}
e^{i\idf(\idt-\idtp)} \quad\quad (\FOT\FOTstar)_{\idf,\idfp} =
\frac{1}{N}\sum_{\idt\in\supt} e^{i\idt(\idf-\idfp)}.
\]
and thus 
\begin{equation}
\label{eq:trG}
\trace(G^*G) = \frac{2(|\supt|\cdot|\supf|)}{N}.
\end{equation}
Observe now that for the null space of $\Phi_\suptf^* \Phi_\suptf $ to
have dimension $K$, at least $K$ of the eigenvalues in \eqref{eq:Gtr2}
must have magnitude greater than or equal to $1$.  As a result
\[
\trace(G^*G) < 2K ~ \Rightarrow ~\dmn(\nullsp(\Phi_\suptf^* \Phi_\suptf )) < K.
\]  
Using the fact that $(a+b)\geq 4ab/(a+b)$ (arithmetic mean dominates
geometric mean), we see that if $|\supt|+|\supf| < 2N$, then
$2|\supt|\cdot|\supf|/N < |\supt|+|\supf|$ which implies
\eqref{eq:Gtr} (and hence $\dmn(\nullsp(\Phi_\suptf^* \Phi_\suptf ))$)
is less than $(|\supt|+|\supf|)/2$.
\end{proof}


\begin{proof}{\bf of Theorem~\ref{th:qrupl0IF}}
  We assume $\suptf$ is selected such that $\Phi_\suptf$ has full
  rank. This happens if $\|\FOTstar \FOT\| < 1$ and
  Theorem~\ref{th:qrupIF} states that this occurs with probability at
  least $1-O(\log N \cdot N^{-\beta})$.
  
  Given this $\suptf$, the (continuous) probability distribution on
  the $\{\alpha(\idtf), \idtf\in\suptf\}$ induces a continuous
  probability distribution on $\range(\Phi_\suptf)$.  We will show
  that for every $\suptfp$ with $|\suptfp|\leq |\suptf|$
\begin{equation}
\label{eq:dims}
\dmn(\range(\Phi_\suptfp)\cap\range(\Phi_\suptf)) < |\suptf|.
\end{equation}
As such, the set of signals in $\range(\Phi_\suptf)$ that have
expansions on a $\suptfp\not=\suptf$ that are {\em at least} as sparse
as their expansions on $\suptf$ is a finite union of subspaces of
dimension strictly smaller than $|\suptf|$.  This set has measure zero
as a subset of $\range(\Phi_\suptf)$, and hence the probability of
observing such a signal is zero.

To show \eqref{eq:dims}, we may assume that $\Phi_\suptfp$ also has full
rank, since if $\dmn(\range(\Phi_\suptfp)) < |\suptfp|$, then
\eqref{eq:dims} is certainly true.  For a set of coefficients $\alpha$
supported on $\suptf$ and $\alphap$ supported on $\suptfp$ to have the
same image under $\Phi$, $\Phi\alpha = \Phi\alphap$ (or equivalently
$\Phi_\suptf\Rtf\alpha = \Phi_\suptfp\Rtfp\alphap$), two things must
be true:
\begin{enumerate}
\item $\alpha$ and $\alphap$ must agree on $\suptf\cap\suptfp$.  This
  is a direct consequence of $\Phi_\suptfp$ being full rank (its
  columns are linearly independent).
\item There is a $\delta\in\nullsp(\Phi)$ such that $\alphap =
  \alpha+\delta$.  Of course,
\[
\delta(\idtf) = 0,\quad \idtf\in(\suptf\cup\suptfp)^c.
\]
By item 1 above, we will also have
\[
\delta(\idtf)=0,\quad \suptf\cap\suptfp.
\]
Thus, $\supp\delta\subset(\suptf\backslash\suptfp)\cup(\suptfp\backslash\suptf)$.
\end{enumerate}
In light of these observations, we see that for
$\dmn(\range(\Phi_\suptfp)\cap\range(\Phi_\suptf)) = |\suptf|$, we
need that for {\em every} $\alpha$ supported on $\suptf$, there is a
$\delta\in\nullsp(\Phi)$ that is supported on
$(\suptf\backslash\suptfp)\cup(\suptfp\backslash\suptf)$ such that 
\[
\delta(\idtf) = -\alpha(\idtf) \quad \idtf\in\suptf\backslash\suptfp.
\]
In other words, we need
\[
\dmn(\nullsp(Q_{(\suptf\backslash\suptfp)\cup(\suptfp\backslash\suptf)}))
\geq \left|\suptf\backslash\suptfp\right|.
\]
However, Lemma~\ref{lm:dimnull} tells us
\begin{eqnarray*}
\dmn(\nullsp(Q_{(\suptf\backslash\suptfp)\cup(\suptfp\backslash\suptf)})) & < &
\frac{\left|\suptf\backslash\suptfp\right|+\left|\suptfp\backslash\suptf\right|}{2} \\
 & \leq & \left|\suptf\backslash\suptfp\right|, 
\end{eqnarray*}
since $|\suptfp|\leq |\suptf|$.  Hence 
$\dmn(\range(\Phi_\suptfp)\cap\range(\Phi_\suptf)) < |\suptf|$, and the theorem follows.
\end{proof}

\subsection{Recovery via $\ell_1$-minimization}
\label{sec:l1}

The problem $(P_0)$ is combinatorial and solving it directly is
infeasible even for modest-sized signals. This is the reason why we
consider instead, the convex relaxation \eqref{eq:(P_1)}.
\begin{theorem}
\label{th:qrupl1IF}
Suppose $f = \Phi \atf$ is a random signal sampled as in
Section~\ref{sec:pmodel} and with parameters obeying 
\begin{equation}
\label{eq:taul1}
\E |T| + \E |\Omega| \le 
\frac{N}{(\beta+1)\log N} \cdot (1/8 + o(1)). 
\end{equation}
Then with probability at least $1- O( (\log N) \cdot N^{-\beta})$, the
solutions of $(P_1)$ and $(P_0)$ are identical and equal to $\atf$.
\end{theorem}
In addition to being computationally tractable, there are analytical
advantages which come with $(P_1)$, as our arguments will essentially
rely on a strong duality result \cite{BVConvex}. In fact, the next
section shows that $\atf$ is a unique minimizer of $(P_1)$ if and only
if there exists a ``dual vector'' $\dv$ satisfying certain properties.
Here, the crucial part of the analysis relies on the fact that
``partial'' Fourier matrices $\FOT := \Rf\FT\Rt^*$ have very
well-behaved eigenvalues, hence the connection with robust uncertainty
principles.

\subsubsection{$\ell_1$-duality}
\label{sec:duality}

For a vector of coefficients $\atf\in\C^{2N}$ supported on
$\supatf:=\supat\cup\supaf$, define the 'sign' vector $\sgn\atf$ by
$(\sgn\atf)(\idtf) := \atf(\idtf)/|\atf(\idtf)|$ for $\idtf\in\supatf$
and $(\sgn\atf)(\idtf) = 0$ otherwise.  We say that $\dv\in C^N$ is a
{\em dual vector} associated to $\atf$ if $\dv$ obeys
\begin{eqnarray}
\label{eq:Psgn}
(\Mtf^*\dv)(\idtf)  =  (\sgn\atf)(\idtf) & \idtf\in\supatf & \\
\label{eq:Plt1}
|(\Mtf^*\dv)(\idtf)| < 1\quad\quad & ~~\idtf\in\supatf^c. &   
\end{eqnarray}
With this notion, we introduce a strong duality result which is
similar to that presented in \cite{CRT}, see also \cite{FuchsDual}.
\begin{lemma}
\label{duality}  
Consider a vector $\alpha\in\C^{2N}$ with support $\supatf =
\supat\cup\supaf$ and put $f = \Mtf\atf$.
\begin{itemize}
\item Suppose that there exists a dual vector and that $\Phi_\Gamma$
  has full rank.  Then the minimizer $\atf^\sharp$ to the problem
  $(P_1)$ is unique and equal to $\atf$.

\item Conversely, if $\atf$ is the unique minimizer of $(P_1)$, then
  there exists a dual vector.
\end{itemize}
\end{lemma}
\begin{proof}
  The program dual to $(P_1)$ is
\begin{equation}
\label{eq:P_1dual}
(D1)\quad\quad
\max_{\dv} ~\Re\left(\dv^*f \right)\quad\quad \mathrm{subject~to}\quad 
\|\Mtf^*\dv\|_{\ell_\infty} \leq 1.
\end{equation}
It is a classical result in convex optimization that if $\tilde{\atf}$
is a minimizer of $(P_1)$, then $\Re(\dv^*\Mtf\tilde{\atf}) \leq
\|\tilde{\atf}\|_{\ell_1}$ for all feasible $\dv$.  Since the primal
is a convex functional subject only to equality constraints, we will
have $\Re(\tilde{\dv}^*\Mtf\tilde{\atf}) = \|\tilde{\atf}\|_{\ell_1}$
if and only if $\tilde{\dv}$ is a maximizer of $(D1)$
\mbox{\cite[Chap. 5]{BVConvex}}.

First, suppose that $\Mtf\Rtf^*$ has full rank and that a dual vector
$\dv$ exists. Set $P = \Mtf^*\dv$.  Then
\begin{eqnarray*}
\Re \<\Mtf\atf,\dv\> & = & \Re\<\atf,\Mtf^*\dv\> \\
 & = & \Re\sum_{\idtf=0}^{N-1} \overline{P(\idtf)}\atf(\idtf) \\
 & = & \Re\sum_{\idtf\in\supatf} \overline{\sgn\atf(\idtf)}\atf(\idtf) \\
 & = & \|\atf\|_{\ell_1}
\end{eqnarray*}
and $\atf$ is a minimizer of $(P_1)$.  Since $|P(\idtf)| < 1$ for
$\idtf\in\supatf^c$, all minimizers of $(P_1)$ must be supported on
$\supatf$.  But $\Mtf\Rtf^*$ has full rank, so $\atf$ is the unique
minimizer.

For the converse, suppose that $\atf$ is the unique minimizer of
$(P_1)$.  Then there exists at least one $\dv$ such that with $P =
\Mtf^*\dv$, $\|P\|_{\ell_\infty} \leq 1$ and $\dv^*f =
\|\atf\|_{\ell_1}$.  Then
\begin{eqnarray*}
\|\atf\|_{\ell_1}  & = & \Re\<\Mtf\atf,\dv\> \\
 & = & \Re\<\atf,\Mtf^*\dv\> \\
 & = & \Re\sum_{\idtf\in\supatf}\overline{P(\idtf)}\atf(\idtf).
\end{eqnarray*}
Since $|P(\idtf)|\leq 1$, equality above can only hold if 
$P(\idtf)=\sgn\atf(\idtf)$ for $\idtf\in\supatf$.

We will argue geometrically that for one of these $\dv$, we have
$|P(\idtf)|<1$ for $\idtf\in\supatf^c$.  Let V be the hyperplane
$\{d\in\C^{2N}:~\Mtf d = f\}$, and let $B$ be the polytope
$B=\{d\in\C^{2N}:~\|d\|_{\ell_1}\leq\|\atf\|_{\ell_1}\}$.  Each of the
$\dv$ above corresponds to a hyperplane $H_\dv =
\{d:~\Re\<d,\Mtf^*\dv\>=\|\atf\|_{\ell_1}\}$ that contains $V$ (since
$\Re\<f,\dv\>=\|\atf\|_{\ell_1}$) and which defines a halfspace
$\{d:~\Re\<d,\Mtf^*\dv\> \leq 1\}$ that contains $B$ (and for each
such hyperplane, a $\dv$ exists that describes it as such).  Since
$\atf$ is the unique minimizer, for one of these $\dvp$, the
hyperplane $H_\dvp$ intersects $B$ only on the minimal facet
$\{d:~\supp d \subset\supatf\}$, and we will have
$P(\idtf)<1,~\idtf\in\supatf^c$.
\end{proof}

Thus to show that $(P_1)$ recovers a representation $\atf$ from a
signal observation $\Phi \alpha$, it is enough to prove that a dual
vector with properties \eqref{eq:Psgn}--\eqref{eq:Plt1} exists.

As a sufficient condition for the equivalence of $(P_0)$ and $(P_1)$,
we construct the {\em minimum energy} dual vector 
\[
\min \|P\|_2, \qquad \text{subject to} \qquad P\in\range(\Phi^*) 
\text{ and } P(\idtf) = \sgn(\alpha)(\idtf), \,\, \forall \idtf \in \suptf.
\]
This minimum energy vector is somehow ``small,'' and we hope that it
obeys the inequality constraints \eqref{eq:Plt1}. Note that $\|P\|_2 =
2 \|S\|_2$, and the problem is thus the same as finding that $S \in
\C^N$ with minimum norm and obeying the constraint above; the solution
is classical and given by
\[
\dv = \Phi_\Gamma (\Phi^*_\Gamma \Phi_\Gamma)^{-1}\Rtf\sgn\alpha
\]   
where again, $\Rtf$ is the restriction operators to $\suptf$.
Setting $P = \Mtf^*\dv$, we need to
establish that
\begin{enumerate}
\item $\Phi^*_\Gamma \Phi_\Gamma$ is invertible (so that $\dv$
  exists), and if so
\item $|P(\idtf)| < 1$ for $\idtf\in\suptf^c$.
\end{enumerate}

The next section shows that for $|\supt|+|\supf|\asymp N/\log N$, not
only is $\Phi^*_\Gamma \Phi_\Gamma$ invertible with high probability
but in addition, the eigenvalues of $(\Phi^*_\Gamma \Phi_\Gamma)^{-1}$
are all less than two, say.  These size estimates will be very useful
to show that $P$ is small componentwise.

\subsubsection{Invertibility}
\label{sec:inverse}

\begin{lemma}
\label{th:inv}
Fix $\beta\geq 1$ and the parameters as in \eqref{eq:taul1}. Then the
matrix $\Phi^*_\Gamma \Phi_\Gamma$ is invertible and obeys
\[
\|(\Phi^*_\Gamma \Phi_\Gamma)^{-1}\| = 1 + o(1).
\]
with probability exceeding $1-O(\log N \cdot N^{-\beta})$.
\end{lemma}
\begin{proof}
  We begin by recalling that with $\FOT$ as before, $\Phi^*_\Gamma
  \Phi_\Gamma$ is given by
\[
\Phi^*_\Gamma \Phi_\Gamma = I + \bpm  0 & 
\FOTstar \\ \FOT & 0 \epm. 
\]
Clearly, $\|(\Phi^*_\Gamma \Phi_\Gamma)^{-1}\| =
1/\lambda_{\mathrm{min}}(\Phi^*_\Gamma \Phi_\Gamma)$ and since
$\lambda_{\mathrm{min}}(\Phi^*_\Gamma \Phi_\Gamma) \ge 1 -
\sqrt{\|\FOTstar\FOT\|}$, we have
\begin{equation*}
\|(\Phi^*_\Gamma \Phi_\Gamma)^{-1}\| \le \frac{1}{1 - \sqrt{\|\FOTstar\FOT\|}}.
\end{equation*}

We then need to prove that $\|\FOTstar\FOT\| = o(1)$ with the required
probability. This follows from the conclusion of Lemma
\ref{lm:auxeigs} which \eqref{eq:taul1} alows to specialize to the
value $1/q = 8\rho_0 \sqrt{(\beta + 1) \log N}$. Note that this gives
more than what is claimed since
\[
\|(\Phi^*_\Gamma \Phi_\Gamma)^{-1}\| \le 1 + \frac{1}{8\rho_0\, 
\sqrt{(\beta + 1) \log N}}  + O(1/\log N). 
\]
\end{proof}

{\bf Remark.} Note that Lemma \ref{lm:auxeigs} assures us that it is
sufficient to take $\E |T| + \E |\Omega|$ of the order of
$N/\sqrt{\log N}$ (rather than of the order of $N/\log N$ as the
Theorem states) and still have invertibility with $\|(\Phi^*_\Gamma
\Phi_\Gamma)^{-1}\| \le 2$, say. The reason why we actually need the
stronger condition will become apparent in the next subsection.

\subsubsection{Proof of Theorem \ref{th:qrupl1IF}}
\label{sec:Pbound}

To prove our theorem, it remains to show that, with high probability,
$|P(\idtf)| < 1$ on $\suptf^c$.
\begin{lemma}
\label{th:Plt1}
Under the hypotheses of Theorem \ref{th:qrupl1IF}, for each
$\idtf\in\suptf^c$
\[
\P\left( |P(\idtf)| \geq 1 \right) ~\leq ~4N^{-(\beta+1)}.
\]
As a result,
\[
\P\left( \max_{\idtf\in\suptf^c} |P(\idtf)| \geq 1\right) ~\leq
~8N^{-\beta}.
\]
\end{lemma}
\begin{proof}
The image of the dual vector $P$ is given by
\begin{equation*}
P := \bpm P_1(\idt) \\ P_2(\idf) \epm
   =  \Phi^* \Phi_\Gamma \, (\Phi^*_\Gamma \Phi_\Gamma)^{-1}\Rtf\sgn\alpha,
\end{equation*}
where the matrix $\Phi^* \Phi_\Gamma$ may be expanded in the time and
frequency subdomains as
\begin{equation*} 
\Phi^* \Phi_\Gamma    =  \bpm \Rt^* & \FT^*\Rf^* \\ \FT\Rt^* & \Rf^* \epm. 
\end{equation*}
Consider first $P_1(\idt)$ for $\idt\in\supt^c$ and let
$V_\idt\in\C^{|\suptf|}$ be the conjugate transpose of the row of the
matrix $\bpm \Rt^* & \FT^*\Rf^* \epm$ corresponding to index $\idt$.
For $\idt\in\supt^c$, the row of $\Rt^*$ with index $\idt$ is zero,
and $V_t$ is then the $(|\supt| + |\supf|)$-dimensional vector
\[
V_\idt = \bpm 0 \\ 
\left\{ \frac{1}{\sqrt{N}}e^{-i\idf\idt},~\idf\in\supf\right\} \epm.
\]
These notations permit to express $P_1(\idt)$ as the inner product
\begin{eqnarray*}
P_1(\idt) & = & 
\<(\Phi^*_\Gamma \Phi_\Gamma)^{-1}\Rtf\sgn\alpha, V_\idt \> \\
& = & \<\Rtf\sgn\alpha, (\Phi^*_\Gamma \Phi_\Gamma)^{-1}V_\idt \> \\
& = & \sum_{\idtf\in\suptf}\overline{W(\idtf)}\sgn\alpha(\idtf)
\end{eqnarray*}
where $W = (\Phi^*_\Gamma \Phi_\Gamma)^{-1}V_\idt$.  The signs of
$\alpha$ on $\suptf$ are statistically independent of $\suptf$ (and
hence of $W$) and, therefore, for a fixed support set $\suptf$,
$P_1(\idt)$ is a weighted sum of independent complex-valued random
variables
\[
P_1(\idt) = \sum_{\idtf\in\suptf} X_\idtf
\]
with $\E X_\idtf = 0$ and $|X_\idtf| \leq |W(\idtf)|$.  Applying the
complex Hoeffding inequality (see the Appendix) gives a
bound on the conditional ditribution of $P(t)$
\[
\P\left( |P_1(\idt)| \geq 1 \,\, | \,\, \suptf \right) ~\leq ~
4\exp\left(-\frac{1}{4\|W\|_2^2} \right).
\]
Thus, it suffices to develop a bound on the magnitude of the vector
$W$. 
Controlling the eigenvalues of  $\|(\Phi^*_\Gamma\Phi_\Gamma)^{-1}\|$ 
is essential here, as
\begin{equation}
\|W\| \le \|(\Phi^*_\Gamma \Phi_\Gamma)^{-1}\| \cdot \|V_\idt\|.
\end{equation}
On the one hand, $\|V_\idt\| = \sqrt{|\supf|/N}$ and as we have seen,
size estimates about $|\supf|$ give $\|V_\idt\| \le
\sqrt{2(p_T+p_\Omega)}$ with the desired probability. On the other
hand, we have also seen that $\|(\Phi^*_\Gamma \Phi_\Gamma)^{-1}\| \le
1 + o(1)$--- also with the desired probability---and,
therefore,
\[
\|W\|^2 \le 2 \cdot \left(1 + o(1))\right) \cdot (p_T +
p_\Omega).
\]
This gives 
\[
P\left( |P_1(\idt)| \geq 1\right) \leq 4\exp\left(-\frac{1}{8 (p_T +
    p_\Omega) (1 + o(1))}\right).
\]
Select $p_T + p_\Omega$ as in \eqref{eq:taul1}. Then
\begin{equation*}
\P\left( |P_1(\idt)| \geq 1\right) 
 \leq  4 \exp(-(\beta+1)\log N)
 \leq  4 N^{-(\beta+1)}
\end{equation*}
and
\[
\P\left( \max_{\idt\in\supt^c}|P_1(\idt)| \geq 1\right) ~\leq ~4N^{-\beta}.
\]
As we alluded earlier, the bound about the size of each individual
$P(t)$ one would obtain assuming that $\E |T| + \E |\Omega|$ be only
of the order $N/\sqrt{\log N}$ would not allow taking the supremum via
the standard union bound. Our approach requires $\E |T| + \E |\Omega|$
to be of the order $N/\log N$.

By the symmetry of the Fourier transform, the same is true for
$P_2(\idf)$. This finishes the proof of Lemma and \ref{th:Plt1} and of
Theorem \ref{th:qrupl1IF}.
\end{proof}

\section{Robust UPs and Basis Pursuit}
\label{sec:BPgeneral}

The results of Sections~\ref{sec:qrup} and \ref{sec:bpIF} extend to
the general situation where the dictionary $\Mtf$ is a union of two
orthonormal bases $\Mt,\Mf$.  In this section, we present results for
pairs of orthogonal bases that parallel those for the time-frequency
dictionary presented in Sections~\ref{sec:qrup} and \ref{sec:bpIF}.
The bounds will depend critically on the degree of similarity of $\Mt$
and $\Mf$, which we measure using the the mutual incoherence defined
in \eqref{eq:mich}, $\ich:=\ich(\Mt,\Mf)$.  As we will see, our
generalization introduces additional ``$\log N$'' factors.  It is our
conjecture that bounds that do not include these factors exist.

As before, the key result is the quantitative robust uncertainty
principle.  We use the same probabilistic setup to sample the support
sets $\supat$, $\supaf$ in the $\Mt$ and $\Mf$ domains
respectively. The statement below is the analogue of
Theorem~\ref{th:qrupIF}. 
\begin{theorem}
\label{th:qrupgen}
Let $\Mtf:=\bpm \Mt & \Mf \epm$ be a dictionary composed of a union of
two orthonormal bases with mutual incoherence $\ich$.  Suppose the
sampling parameters obey
\begin{equation}
\label{eq:taugenup}
\E |\Gamma_1| + \E |\Gamma_2| 
\leq \frac{C_1}{\ich^2\cdot((\beta+1)\log N)^{5/2}}
\end{equation}
for some positive constant $C_1 > 0$.  Assume $\mu \le 1/\sqrt{2(\beta
  +1)\log N}$.  Then with probability at least $1-O(\log N \cdot
N^{-\beta})$, every signal $f$ with $\Mt f$ supported on $\supat$ has
most of its energy in the $\Mf$-domain outside of $\supaf$:
\[
\|\Mf f \cdot 1_{\Gamma_2}\|^2 \leq \|f\|^2/2,
\] 
and vice versa.  As a result, for nearly all pairs $(\Mt,\Mf)$ with
sizes obeying \eqref{eq:taugenup}, it is impossible to find a signal
$f$ supported on $\supat$ in the $\Mt$-domain and $\supaf$ in the
$\Mf$-domain.
\end{theorem}
We would like to re-emphasize the significant difference between these
results and \eqref{eq:gendiscreteup}. Namely, \eqref{eq:taugenup}
effectively squares the size of the joint support since, ignoring
log-like factors, the factor $1/\mu$ is replaced by $1/\mu^2$.  For
example, in the case where the two bases are maximally incoherent,
i.e. $\mu = 1/\sqrt{N}$, our condition says that it is nearly
impossible to concentrate a function in both domains simultaneously
unless (again, up to logarithmic factors)
\[
|\supat| + |\supaf| \sim N, 
\]
which needs to be compared with \eqref{eq:gendiscreteup}
\[
|\supat| + |\supaf| \geq 2 \sqrt{N}. 
\]
For mutual incoherences scaling like a power-law $\mu \sim
N^{-\gamma}$, our condition essentially reads $|\supat| + |\supaf|
\sim N^{2\gamma}$ compared to $|\supat| + |\supaf| \sim N^{\gamma}$.

The proof of Theorem~\ref{th:qrupgen} directly parallels that of
Theorem~\ref{th:qrupIF}, with $A:= \Raf\Mf^*\Mt\Rat^*$ playing the
role of the partial Fourier transform from $T$ to $\Omega$. Our
argument calls for bounds on the eigenvalues of the random matrix $A^*
A$ which we write as the sum of two terms; a diagonal and an
off-diagonal term
\[
A^* A = D + {\cal H}_1. 
\]
We use large deviation theory to control the norm of $D$ while bounds
on the size of ${\cal H}_1$ are obtained by using moment estimates.
This calculation involves estimates about the expected value of the
Frobenius norm of large powers of $A^* A$ and is very delicate. We do
not reproduce all these arguments here (this is the scope of a whole
separate article) and simply state a result which is proved in
\cite{SparsityIncoherence}
\begin{equation}
  \label{eq:A*A}
\P(\|A^* A\| \ge 1/2) \le C \cdot \log N \cdot N^{-\beta}
\end{equation}
for $\E |\Gamma_1| + \E |\Gamma_2|$ obeying \eqref{eq:taugenup} (here
$C$ is some universal positive constant). Now for \eqref{eq:A*A} to
hold, we also need that the incoherence be not too large and obeys
$\mu > 1/\sqrt{2(\beta +1)\log N}$ which is the additional condition
stated in the hypothesis. The idea that $\mu$ cannot be too large is
somewhat natural as otherwise for $\mu = 1$, say, the two bases would
share at least one element and we would have $\|A^*A\| = 1$ as soon as
$\Gamma_1$ and $\Gamma_2$ would contain a common element.  As we have
seen in section 3, the size estimate \eqref{eq:A*A} would then
establish the theorem.

The generalized $\ell_0$-uniqueness result
follows directly from Theorem~\ref{th:qrupgen}:
\begin{theorem}
\label{th:qrupl0gen}
Let $f = \Phi \alpha$ be an observed signal sampled as in
Section~\ref{sec:pmodel}, and with parameters obeying 
\[
\E |\Gamma_1| + \E |\Gamma_2| \le \frac{C_2}{\ich^2\cdot((\beta+1)\log
  N)^{5/2}}.
\]
Assume $\mu \le 1/\sqrt{2(\beta +1)\log N}$. Then with probability
$1-O(\log N \cdot N^{-\beta})$, the solution to $(P_0)$ is unique and
equal to $\atf$.
\end{theorem}
The only change to the proof presented in Section~\ref{sec:l0uniq} is
in the analogue to Lemma~\ref{lm:dimnull}:
\begin{lemma}
\label{lm:dimnullgen}
Let $\supat,\supaf$ be fixed subsets of $\{0,\ldots,N-1\}$, let
$\supatf=\supat\cup\supaf$, and let $Q_\supatf$ be the $N\times
|\supatf|$ matrix
\[
Q_\supatf = \bpm \Mt\Rat^* & \Mf\Raf^* \epm.
\]
If $|\supatf| < 2/\ich^2$, then
\[
\dmn(\nullsp(Q_\supatf)) < \frac{|\supatf|}{2}.
\]
\end{lemma}
The proof of Lemma~\ref{lm:dimnullgen} has exactly the same structure
as the proof of Lemma~\ref{lm:dimnull}.  The only modification comes in
calculating the trace of $G^*G$; here each term can be bounded by
$\ich^2$, and we have $\trace(G^*G) \leq 2(|\supat|\cdot
|\supaf|)\ich^2$.  Lemma~\ref{lm:dimnullgen} follows.

The conditions for the equivalence of $(P_0)$ and $(P_1)$ can also be generalized.
\begin{theorem}
\label{th:qrupl1gen}
Let $f = \Mtf\atf$ be a random signal generated as in
Section~\ref{sec:pmodel} with
\[
\E |\Gamma_1| + \E |\Gamma_2| \le \frac{C_3}{\ich^2\cdot
  ((\beta+1)\log N)^{5/2}}.
\] 
Assume $\mu \le 1/\sqrt{2(\beta +1)\log N}$. Then with probability
$1-O(\log N \cdot N^{-\beta})$, the solutions of $(P_0)$ and $(P_1)$
are identical and equal to $\atf$.
\end{theorem}
The proof of Theorem~\ref{th:qrupl1gen} is again almost exactly the
same as that we have already seen.  Using Theorem~\ref{th:qrupgen},
the eigenvalues of $(\Phi^*_\Gamma \Phi_\Gamma)^{-1}$ are controlled,
allowing us to construct a dual vector meeting the conditions
\eqref{eq:Psgn} and \eqref{eq:Plt1} of Section~\ref{sec:duality}.
Note that the $(\log N)^{5/2}$ term in the denominator means that that
$\P(|P(\idtf)| < 1),~\idtf\in\supatf^c$ goes to zero at a much faster
speed than a negative power of $N$, it decays as $\exp(- \rho (\log
N)^5)$ for some positive constant $\rho > 0$.


\section{Numerical Experiments}
\label{sec:numerical}

From a practical standpoint, the ability of $(P_1)$ to recover sparse
decompositions is nothing short of amazing. To illustrate this fact,
we consider a $256$ point signal composed of $60$ spikes and $60$
sinusoids; $|\supt|+|\supf|\approx N/2$, see Figure 1.  Solving
$(P_1)$ recovers the original decomposition {\em exactly}.

We then empirically validate the previous numerical result by
repeating the experiment for various signals and sample sizes, see
Figure~\ref{fig:reccurves}. These experiments were designed as
follows:
\begin{enumerate}
\item set $N_\suptf$ as a percentage of the signal length $N$;
\item select a support set $\suptf = \supt\cup\supf$ of size
  $|\suptf|=N_\suptf$ uniformly at random;
\item sample a vector $\atf$ on $\suptf$ with independent and 
identically distributed
  Gaussian entries\footnote{The results presented here do not seem to depend
    on the actual distribution used to sample the coefficients.};
\item make $f=\Mtf\atf$;
\item solve $(P_1)$ and obtain $\hat{\atf}$;
\item compare $\atf$ to $\hat{\atf}$;
\item repeat $100$ times for each $N_\suptf$;
\item repeat for signal lengths $N=256,512,1024$. 
\end{enumerate}
Figure~\ref{fig:reccurves}(a) shows that we are numerically able to
recover ``sparse'' superpositions of spikes and sinusoids when
$|\supt|+|\supf|$ is close to $N/2$, at least for this range of sample sizes
$N$ (we use the quotations since decompositions of this order can
hardly be considered sparse).  Figure~\ref{fig:reccurves}(b) plots the
success rate of the sufficient condition for the recovery of the
sparsest $\atf$ developed in Section~\ref{sec:duality} (i.e. the
minimum energy signal is a dual vector).  Numerically, the sufficient
condition holds when $|\supt|+|\supf|$ is close to $N/5$.

The time-frequency dictionary is special in that it is maximally
incoherent ($\ich = 1$).  But as suggested in \cite{DonohoHuo},
incoherence between two bases is the rule, rather than the exception.
To illustrate this, the above experiment was repeated for $N=256$ with
a dictionary that is a union of the spike basis and of an orthobasis
sampled uniformly at random (think about orthogonalizing $N$ vectors
sampled independently and uniformly on the unit-sphere of $\C^N$).  As
shown in Figure~\ref{fig:reccurvesrandom}, the results are very close
to those obtained with time-frequency dictionaries; we recover
``sparse'' decompositions of size about $|\supat| + |\supaf| \leq 0.4\cdot
N$.

\newcommand{\Figures}{Figures}

\begin{figure}
\label{fig:recex}
\centerline{
\begin{tabular}{cccccc}
$\atf$ & $f = \Mtf\alpha$ & & spike component & & sinusoidal component \\
\includegraphics[width=1.5in]{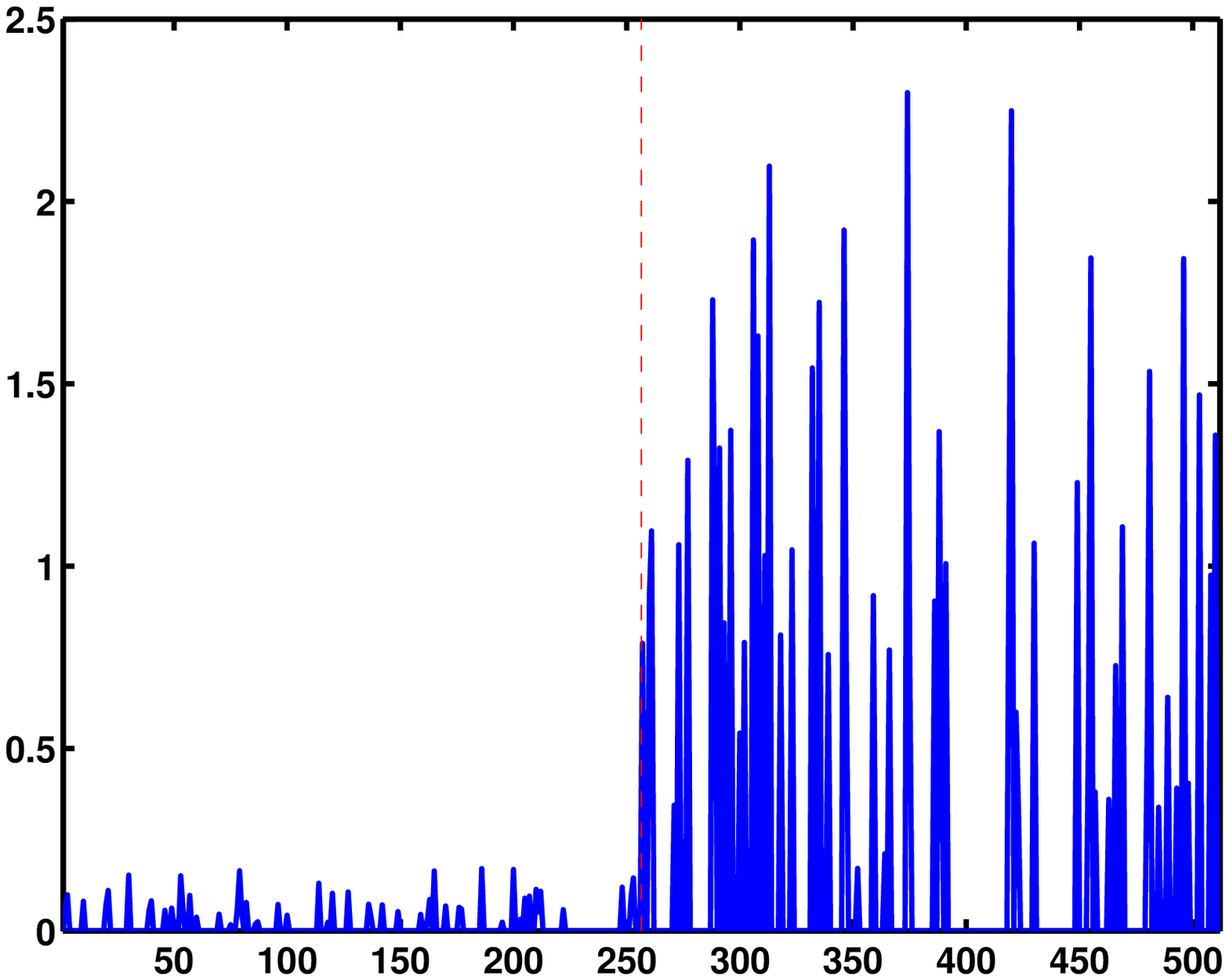} &
\includegraphics[width=1.5in]{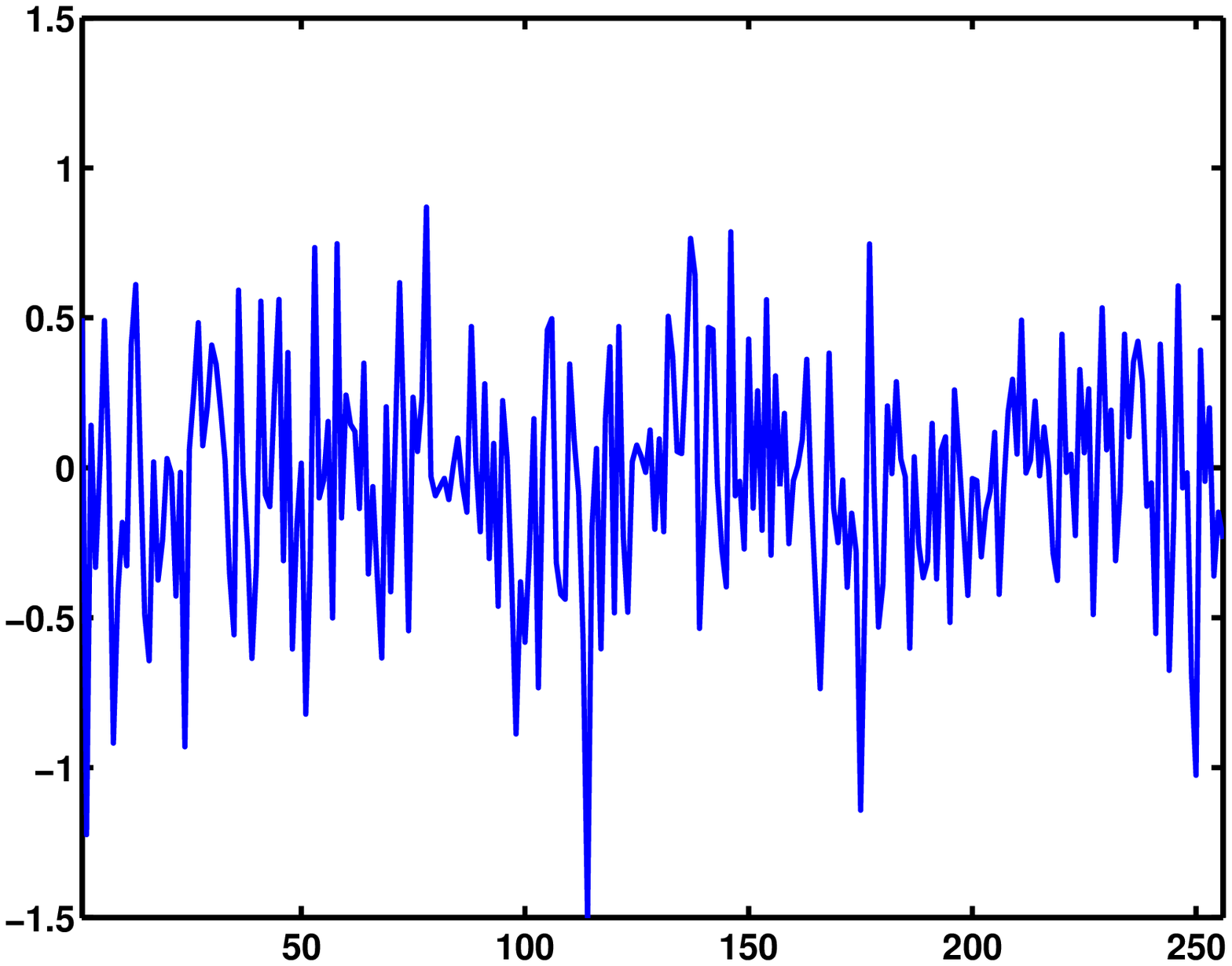} &
\raisebox{0.6in}{{\bf =}} &
\includegraphics[width=1.5in]{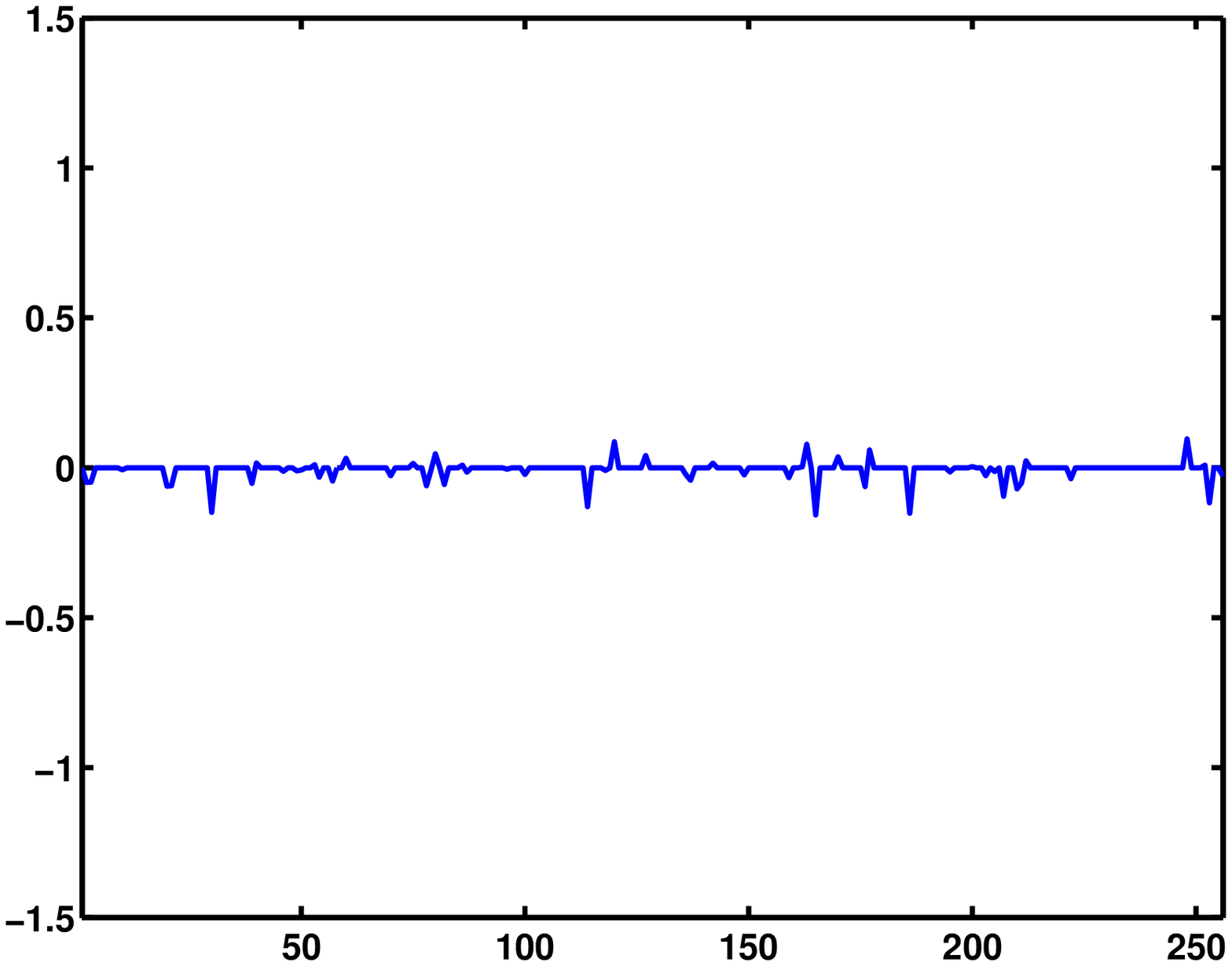} &
\raisebox{0.6in}{{\bf +}} &
\includegraphics[width=1.5in]{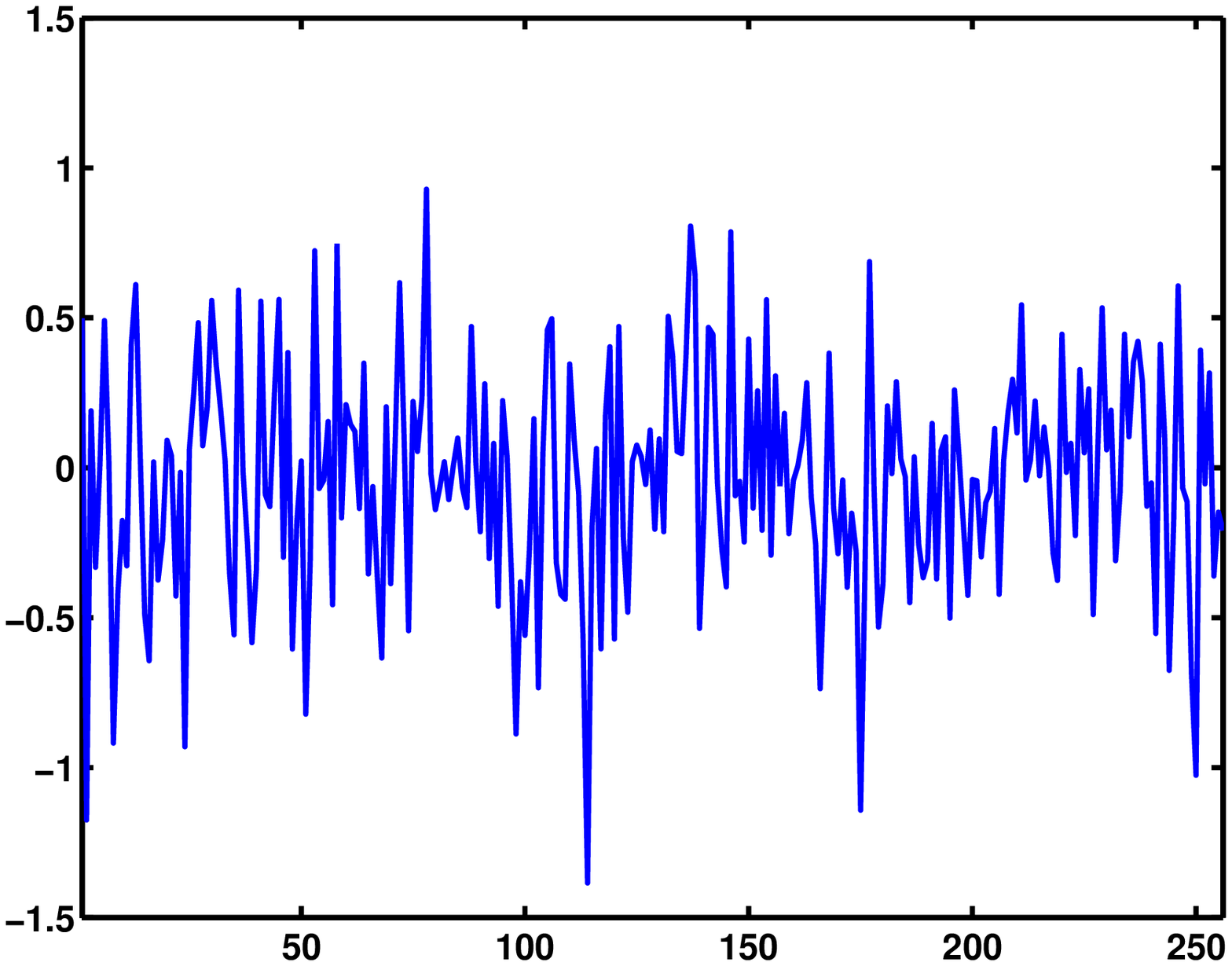} \\
(a) & (b) & & (c) & & (d) 
\end{tabular}
}
\caption{\small\sl Recovery of a ``sparse'' decomposition. 
  (a) Magnitudes of a randomly generated coefficient vector $\atf$ with
  $120$ nonzero components.  The spike components are on the left
  (indices 1--256) and the sinusoids are on the right (indices
  257--512).  The spike magnitudes are made small compared to the
  magnitudes of the sinusoids for effect; we cannot locate the spikes
  by inspection from the observed signal $f$, whose real part is shown
  in (b).  Solving $(P_1)$ separates $f$ into its spike (c) and
  sinusoidal components (d) (the real parts are plotted).  }
\end{figure}

\begin{figure}
\centerline{
\begin{tabular}{ccc}
$\ell_1$-recovery & & sufficient condition \\
\raisebox{0.8in}{\rotatebox{90}{\% success}}
\includegraphics[width=3in]{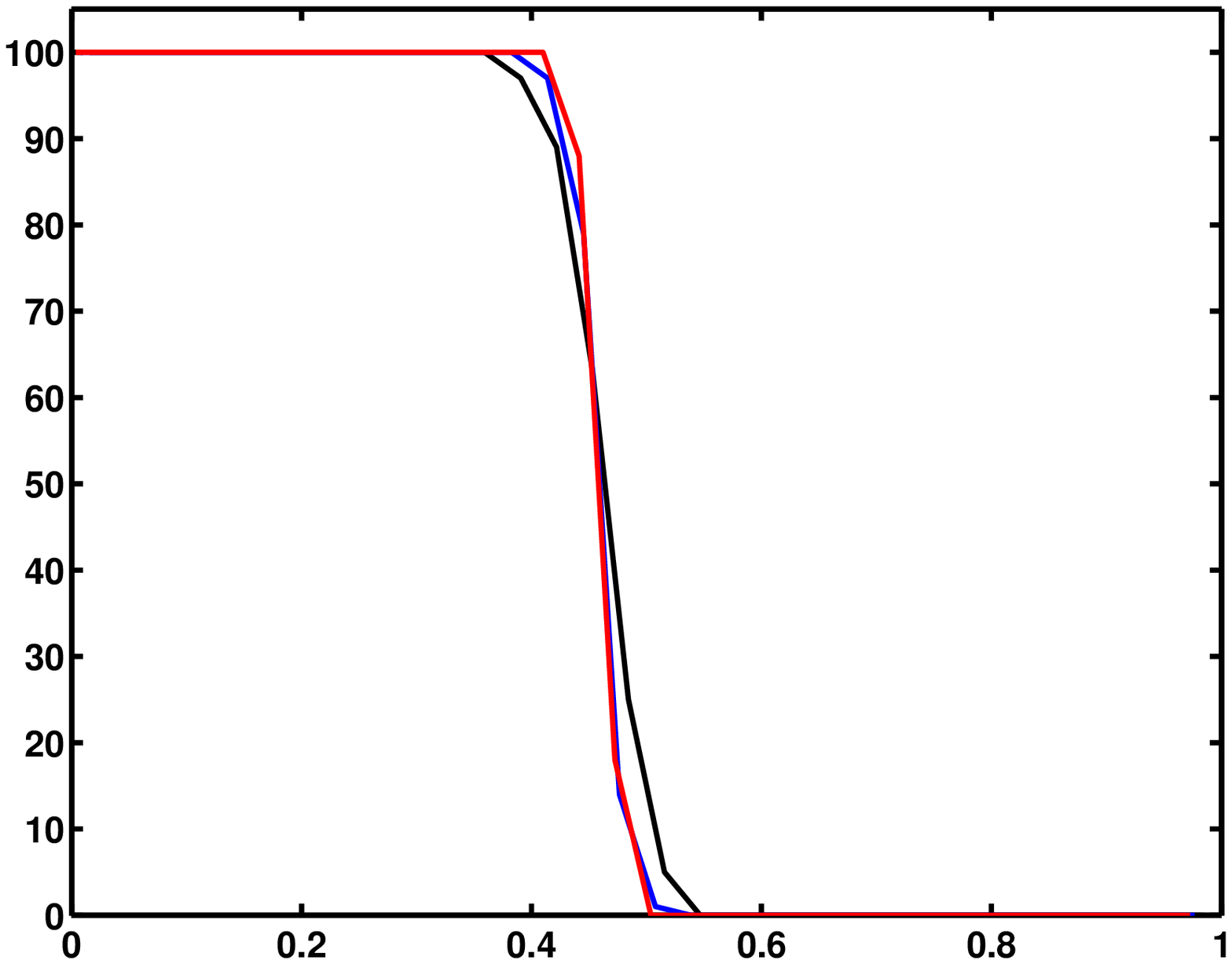} & \hspace{5mm} &
\raisebox{0.8in}{\rotatebox{90}{\% success}}
\includegraphics[width=3in]{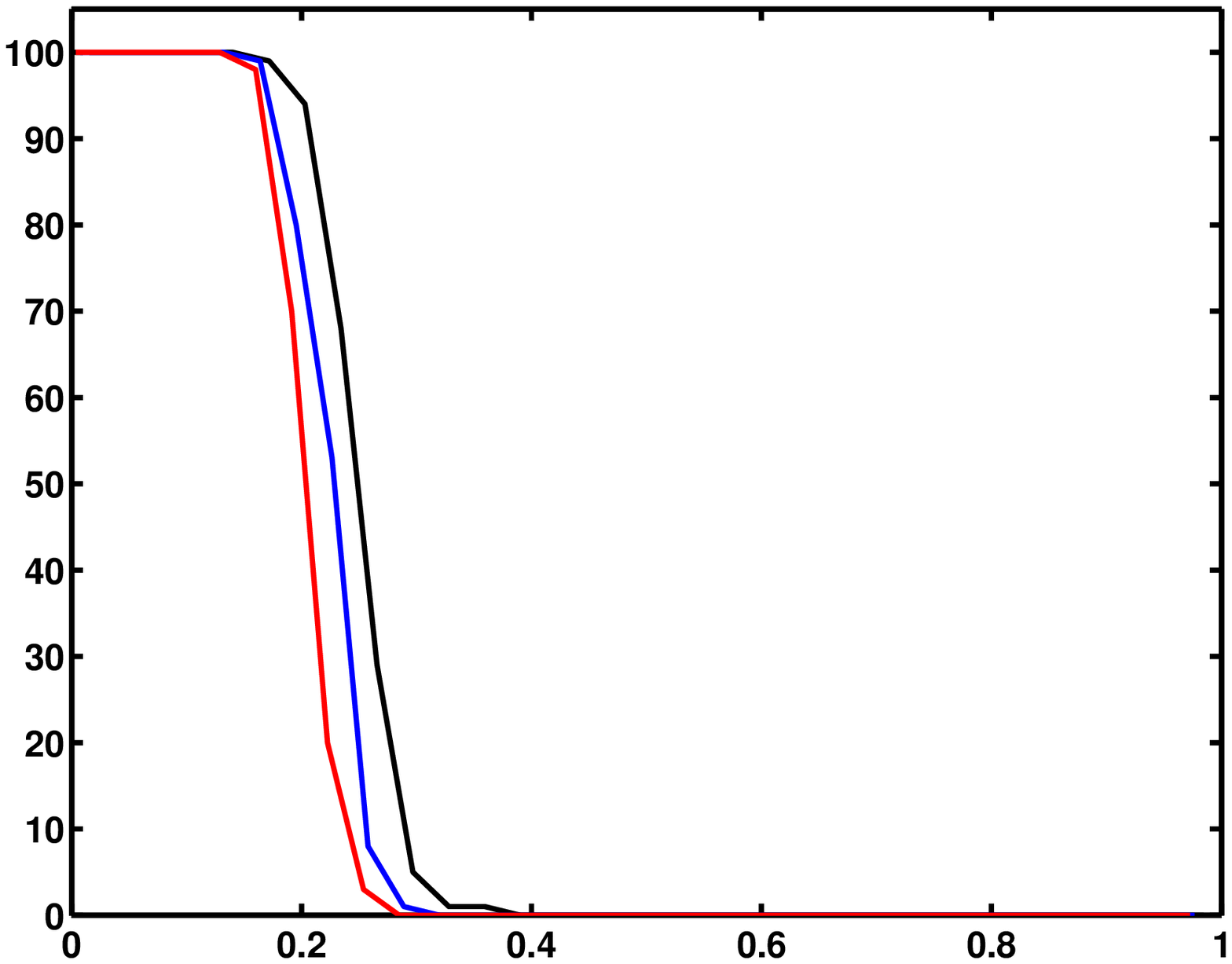} \\
$(|\supt|+|\supf|)/N$ & & $(|\supt|+|\supf|)/N$
\end{tabular}
}
\caption{\small\sl $\ell_1$-recovery for the time-frequency dictionary.  
  (a) Success rate of $(P_1)$ in recovering the sparsest decomposition
  versus the number of nonzero terms.  (b) Success rate of the
  sufficient condition (the minimum energy signal is a dual vector).
}
\label{fig:reccurves}
\end{figure}

\begin{figure}
\centerline{
\begin{tabular}{ccc}
$\ell_1$ recovery & & sufficient condition \\
\raisebox{0.8in}{\rotatebox{90}{\% success}}
\includegraphics[width=3in]{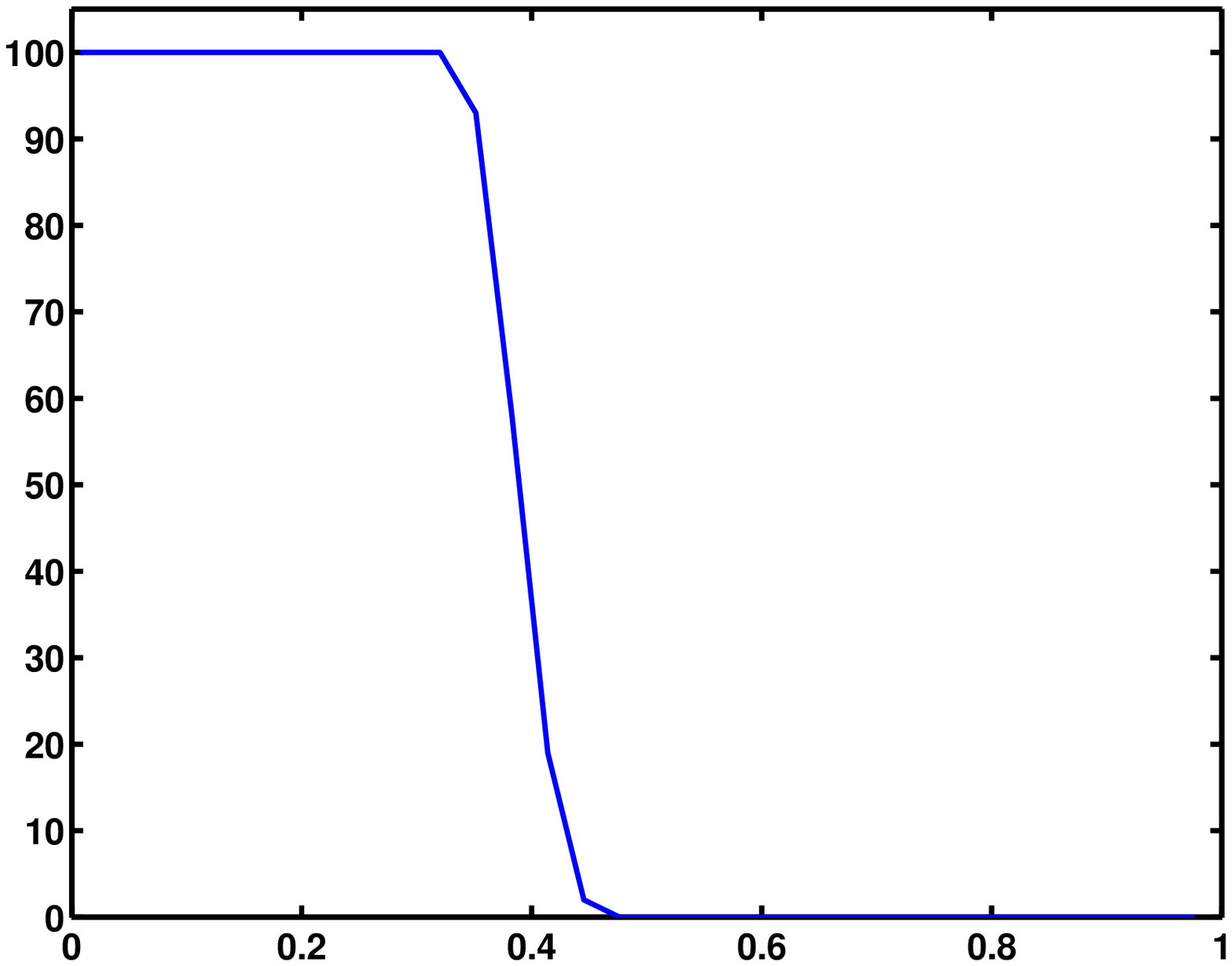} & 
\hspace{5mm} &
\raisebox{0.8in}{\rotatebox{90}{\% success}}
\includegraphics[width=3in]{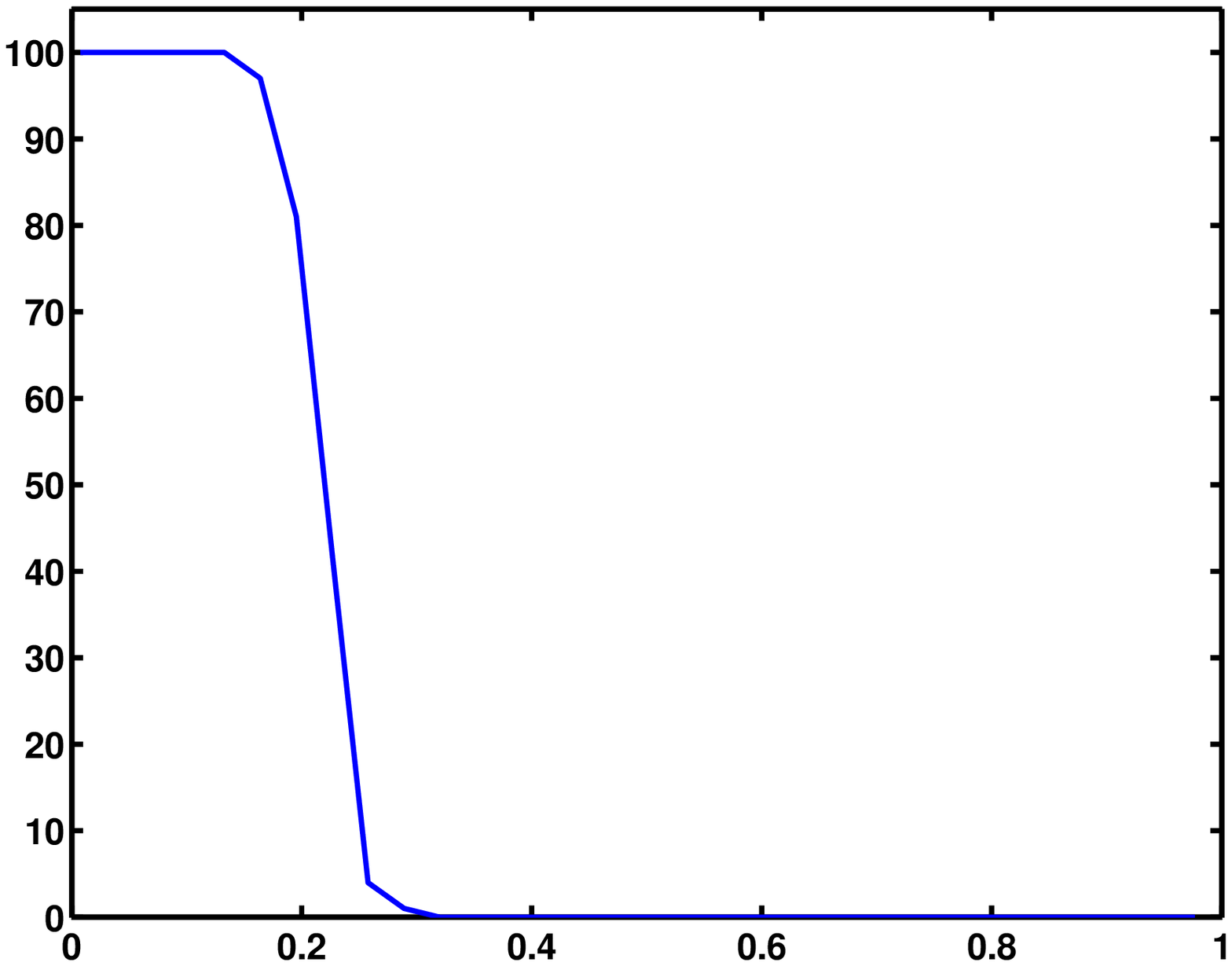} \\
$(|\supt|+|\supf|)/N$ & & $(|\supt|+|\supf|)/N$
\end{tabular}
}
\caption{\small\sl $\ell_1$-recovery for the spike-random dictionary.  
  (a) Success rate of $(P_1)$ in recovering the sparsest decomposition
  versus the number of nonzero terms.  (b) Success rate of the
  sufficient condition.  }
\label{fig:reccurvesrandom}
\end{figure}

\section{Discussion}
\label{sec:discussion}

In this paper, we have demonstrated that except for a negligible
fraction of pairs $(T, \Omega)$, the behavior of the discrete
uncertainty relation is very different from what worst case
scenarios---which have been the focus of the literature thus
far---suggest.  We introduced probability models and a robust
uncertainty principle showing that for for nearly all pairs $(T,
\Omega)$, it is actually impossible to concentrate a discrete signal
on $T$ and $\Omega$ simultaneously unless the size of the joint
support $|T| + |\Omega|$ be at least of the order of $N/\sqrt{\log
  N}$.  We derived significant consequences of this new uncertainty
principle, showing how one can recover sparse decompositions by
solving simple convex programs.

Our sampling models were selected in perhaps the most natural way,
giving to each time point and to each frequency point the same chance
of being sampled, independently of the others. Now there is little
doubt that conclusions similar to those derived in this paper would
hold for other probability models. In fact, our analysis develops a
machinery amenable to other setups. The centerpiece is the study of
the singular values of partial Fourier transforms. For other sampling
models such as models biased toward low or high ferequencies for
example, one would need to develop analogues of Lemma
\ref{lm:auxeigs}.  Our machinery would then nearly automatically
transforms these new estimates into corresponding claims.

In conclusion, we would like to mention areas for possible improvement
and refinement.  First, although we have made an effort to obtain
explicit constants in all our statements (with the exception of
section \ref{sec:BPgeneral}), there is little doubt that a much more
sophisticated analysis would yield better estimates for the singular
values of partial Fourier transforms, and thus provide better
constants. Another important question we shall leave for future
research, is whether the $1/\sqrt{\log N}$ factor in the QRUP (Theorem
\ref{th:qrupIF}) and the $1/\log N$ for the exact
$\ell_1$-reconstruction (Theorem \ref{th:qrupl1IF}) are necessary.
Finally, we already argued that one really needs to randomly sample
the support to derive our results but we wonder whether one needs to
assume that the signs of the coefficients $\alpha$ (in $f = \Phi
\alpha$) need to be randomized as well.  Or would it be possible to
show analogs of Theorem \ref{th:qrupl1IF} ($\ell_1$ recovers the
sparsest decomposition) for all $\alpha$, provided that the support of
$\alpha$ may not be too large (and randomly selected)?  Recent work
\cite{OptimalRecovery,LPdecode} suggests that this might be
possible---at the expense of additional logarithmic factors.

\section{Appendix: Concentration-of-Measure Inequalities}

The Hoeffding inequality is a well-known large deviation bound for
sums of independent random variables.  For a proof and interesting
discussion, see \cite{Lugosi-Notes}.

\begin{lemma}
\label{lm:hoeffding}
(Hoeffding inequality) Let $X_0,\ldots,X_{N-1}$ be independent
real-valued random variables such that $\E X_j = 0$ and $|X_j| \leq
a_j$ for some positive real numbers $a_j$.  For $\epsilon > 0$
\[
\P\left(\left| \sum_{j=0}^{N-1} X_j \right| \geq \epsilon\right) \leq
2\exp\left(-\frac{\epsilon^2}{2\|\mathbf{a}\|^2_2}\right) 
\]
where $\|\mathbf{a}\|^2_2 = \sum_j a^2_j$.
\end{lemma}

\begin{lemma}
\label{lm:choeffding}
(complex Hoeffding inequality) Let $X_0,\ldots,X_{N-1}$ be independent
complex-valued random variables such that $\E X_j = 0$ and $|X_j| \leq
a_j$.  Then for $\epsilon > 0$
\[
\P\left(\left|\sum_{j=0}^{N-1} X_j \right|\geq \epsilon\right) \leq
4\exp\left(-\frac{\epsilon^2}{4\|\mathbf{a}\|^2_2}\right).
\]
\end{lemma}
\begin{proof}
  Separate the $X_j$ into their real and imaginary parts; $X^\re_j =
  \Re X_j,~X^\im_j = \Im X_j$.  Clearly, $|X^\re_j| \leq a_j$ and
  $|X^\im_j| \leq a_j$.  The result follows immediately from
  Lemma~\ref{lm:hoeffding} and the fact that
\[
\P\left(\left|\sum_{j=0}^{N-1}X_j \right|\geq\epsilon\right) \leq
\P\left(\left|\sum_{j=0}^{N-1}X_j^\re \right|\geq\epsilon/\sqrt{2}\right) + 
\P\left(\left|\sum_{j=0}^{N-1}X_j^\im \right|\geq\epsilon/\sqrt{2}\right).
\]
\end{proof}


\begin{thebibliography}{11}
  
\bibitem{MassartSharp} S. Boucheron, G. Lugosi, and P. Massart, A
  sharp concentration inequality with applications, {\em Random
  Structures Algorithms} {\bf 16} (2000), 277--292.

\bibitem{BVConvex} S. Boyd, and L. Vandenberghe, {\em Convex
    Optimization}, Cambridge University Press, 2004.

\bibitem{CRT} E. J. Cand\`es, J. Romberg, and T. Tao, Robust
  uncertainty principles: exact signal reconstruction from highly
  incomplete frequency information, Technical Report, California
  Institute of Technology.  Submitted to {\em IEEE Transactions on
    Information Theory}, June 2004.  Available on the ArXiV preprint server: {\tt 
    math.GM/0409186}. 

  
\bibitem{SparsityIncoherence} E. J. Cand\`es, and J. Romberg, The Role of
  Sparsity and Incoherence for Exactly Reconstructing a Signal from
  Limited Measurements, Technical Report, California Institute of
  Technology.
      
\bibitem{OptimalRecovery} E. J. Cand\`es, and T. Tao, Near optimal
  signal recovery from random projections: universal encoding
  strategies?  Submitted to {\em IEEE Transactions on Information
    Theory}, October 2004. Available on the ArXiV preprint server: {\tt 
    math.CA/0410542}. 
  
\bibitem{LPdecode} E. J. Cand\`es, and T. Tao, Decoding of random
  linear codes.  Manuscript, October 2004. 

\bibitem{BP} S. S. Chen, D. L. Donoho, and M. A. Saunders, Atomic
  decomposition by basis pursuit, {\em SIAM J. Scientific Computing}
  \textbf{20} (1999), 33--61.
  
\bibitem{DonohoGeomSep} D. L. Donoho, {\em Geometric separation using
  combined curvelet/wavelet representations}. Lecture at the
  International Conference on Computational Harmonic Analysis,
  Nashville, Tennessee, May 2004. 
  
\bibitem{DonohoStark} D. L. Donoho, P. B. Stark, Uncertainty principles
  and signal recovery, {\em SIAM J. Appl. Math.} \textbf{49} (1989),
  906--931.
  
\bibitem{DonohoHuo} D. L. Donoho and X. Huo, Uncertainty principles and
  ideal atomic decomposition, {\em IEEE Transactions on Information
    Theory}, \textbf{47} (2001), 2845--2862.
  
\bibitem{DonohoElad} D. L. Donoho and M. Elad, Optimally sparse
  representation in general (nonorthogonal) dictionaries via $\ell_1$
  minimization.  {\em Proc. Natl. Acad. Sci. USA} \textbf{100} (2003),
  2197--2202. 
  
\bibitem{EladBruckstein} M. Elad and A. M. Bruckstein, A generalized
  uncertainty principle and sparse representation in pairs of $\R^N$
  bases, {\em IEEE Transactions on Information Theory}, \textbf{48}
  (2002), 2558--2567.
  
\bibitem{FeuerNemirovsky} A. Feuer and A. Nemirovsky, On sparse
  representations in pairs of bases, Accepted to the {\em IEEE
    Transactions on Information Theory} in November 2002.

\bibitem{FuchsDual} J. J. Fuchs, On sparse representations in arbitrary
  redundant bases, {\em IEEE Transactions on Information Theory},
  \textbf{50} (2004), 1341--1344.
  
\bibitem{GribonvalNielsen} R. Gribonval and M. Nielsen, Sparse
  representations in unions of bases. {\em IEEE Trans. Inform. Theory}
  {\bf 49} (2003), 3320--3325. 
        
\bibitem{prolateII} H. J. Landau and H. O. Pollack, Prolate spheroidal
  wave functions, Fourier analysis and uncertainty II, {\em Bell
    Systems Tech. Journal}, \textbf{40} (1961), pp. 65--84.
  
\bibitem{prolateIII} H. J. Landau and H. O. Pollack, Prolate spheroidal
  wave functions, Fourier analysis and uncertainty III, {\em Bell
    Systems Tech. Journal}, \textbf{41} (1962), pp. 1295--1336.
    
\bibitem{Lugosi-Notes} G. Lugosi, Concentration-of-measure
  Inequalities, Lecture Notes.

\bibitem{MeyerAverbuchCoifman} F. G. Meyer, A. Z. Averbuch and R. R.
  Coifman, Multi-layered Image Representation: Application to
    Image Compression, {\em IEEE Transactions on Image
    Processing}, {\bf 11} (2002), 1072-1080. 
    
\bibitem{prolateI} D. Slepian and H. O. Pollak, Prolate spheroidal wave
  functions, Fourier analysis and uncertainty I, {\em Bell Systems
    Tech. Journal}, \textbf{40} (1961), pp. 43--64.
  
\bibitem{prolateV} D. Slepian, Prolate spheroidal wave functions,
  Fourier analysis and uncertainty V --- the discete case, {\em Bell
    Systems Tech. Journal}, \textbf{57} (1978), pp. 1371--1430.
  
\bibitem{StarckAstro} J. L.  Starck, E. J. Cand\`es, and D. L. Donoho.
  Astronomical image representation by the curvelet transform {\em
    Astronomy \& Astrophysics}, {\bf 398} 785--800.
  
\bibitem{EdgesTextures} J. L. Starck, M. Elad, and D. L. Donoho.  Image
  Decomposition Via the Combination of Sparse Representations and a
  Variational Approach. Submitted to {\em IEEE Transactions on Image
    Processing}, 2004.

\bibitem{stev-lenstra} P. Stevenhagen, H. W. Lenstra Jr., Chebotar\"ev
  and his density theorem, {\em Math. Intelligencer} \textbf{18}
  (1996), no. 2, 26--37.
  
\bibitem{tao:uncertainty} T. Tao, An uncertainty principle for cyclic
  groups of prime order, preprint. {\rm math.CA/0308286}

\bibitem{Tropp03} J. A. Tropp, Greed is good: {A}lgorithmic results for sparse
  approximation, Technical Report, The University of Texas at Austin, 2003.


\end{thebibliography}
\end{document}